\documentstyle{article}[12pt]

\def\be{\begin{equation}}
\def\ee{\end{equation}}

\def\a{\alpha}
\def\b{\beta}
\def\d{\delta}
\def\G{\Gamma}
\def\lb{\lambda}
\def\s{\sigma}
\def\t{\theta}
\def\ve{\varepsilon}
\def\vp{\varphi}

\def\Ac{{\cal A}}
\def\Bc{{\cal B}}
\def\Hc{{\cal H}}

\def\Jc{{\cal J}}
\def\\Kc{{\cal K}}
\def\Lc{{\cal L}}
\def\Sc{{\cal S}}
\def\Zc{{\zeta}}

\def\fl{\forall \, }
\def\ify{\infty}
\def\ov{\overline}
\def\longra{\longrightarrow}
\def\sbs{\subset}
\def\ts{\times}
\def\wdg{\wedge}
\def\wh{\widehat}

\def\build#1_#2^#3{\mathrel{
\mathop{\kern 0pt#1}\limits_{#2}^{#3}}}

\font\tenbb=msbm10
\font\sevenbb=msbm7
\font\fivebb=msbm5
\newfam\bbfam
\textfont\bbfam=\tenbb 
\scriptfont\bbfam=\sevenbb
\scriptscriptfont\bbfam=\fivebb
\def\bb{\fam\bbfam}

\def\Nb{{\bb N}}
\def\Tb{{\bb T}}
\def\Zb{{\bb Z}}
\def\Cb{{\bb C}}
\def\Rb{{\bb R}}

\def\wh{\widehat}

\def\Ac{{\cal A}}
\def\Bc{{\cal B}}
\def\Cc{{\cal C}}
\def\Dc{{\cal D}}
\def\Hc{{\cal H}}
\def\Lc{{\cal L}}

\def\Sc{{\cal S}}
\def\Uc{{\cal U}}

\def\Hc{{\cal H}}

\def\Kc{{\cal K}}
\def\Lc{{\cal L}}

\def\Sc{{\cal S}}
\def\Uc{{\cal U}}
\newcommand{\wt}{\widetilde}

\def\a{\alpha}
\def\b{\beta}
\def\d{\delta}
\def\lb{\lambda}
\def\g{\gamma}
\def\om{\omega}
\def\s{\sigma}
\def\t{\theta}
\def\ve{\varepsilon}
\def\vp{\varphi}

\def\z{\zeta}

\def\D{\Delta}
\def\G{\Gamma}
\def\Lb{\Lambda}
\def\Om{\Omega}
\def\Si{\Sigma}

\def\fl{\forall \, }
\def\ify{\infty}
\def\op{\oplus}
\def\ot{\otimes}
\def\ov{\overline}
\def\ra{\rightarrow}
\def\longra{\longrightarrow}
\def\sbs{\subset}
\def\ts{\times}
\def\wdg{\wedge}
\def\wh{\widehat}

\def\fl{\forall}
\def\ify{\infty}

\def\nb{\nabla}
\def\op{\oplus}
\def\ot{\otimes}
\def\ov{\overline}
\def\part{\partial}

\def\sbs{\subset}
\def\semi{>\!\!\!\lhd}

\def\ts{\times}
\def\wdg{\wedge}
 \title{{\bf  Cyclic Cohomology, Quantum group Symmetries and the Local Index Formula for} $SU_q(2)$
}
\author{Alain Connes }
 
\date{}

\begin{document}

\maketitle

\begin{abstract}
\noindent We analyse the NC-space
underlying the quantum group
$SU_q(2)$ from the spectral point of view which
is the basis of noncommutative geometry,
and show how the general theory developped in our 
joint work with H. Moscovici 
applies to the specific spectral triple 
defined by Chakraborty and Pal. 
This provides the pseudo-differential calculus, the 
 Wodzciki-type residue, and the local cyclic cocycle giving
the index formula. The cochain whose coboundary
is the difference between the original 
Chern character and the local one is given by the remainders in 
the rational approximation of the logarithmic derivative of the 
Dedekind eta function.
This specific example allows to illustrate
the general notion of locality in NCG. The formulas computing the residue 
are "local". Locality
 by stripping all the expressions from
irrelevant details makes them computable.
The key feature of this spectral triple
is its equivariance, i.e. the $SU_q(2)$-symmetry.
We shall explain how this leads naturally
to the general concept of invariant cyclic 
cohomology in the framework of quantum group symmetries. 

\end{abstract}
 %%%%\begin{center}
%%%%{ Contents}
%%%%\end{center}

\begin{enumerate}
\item[1.] {\bf Introduction}
\item[2.]{\bf  Operator theoretic Local Index Formula}
\item[3.] {\bf Dimension Spectrum of ${\rm SU}_q (2)$, ($q=0$)}
\item[4.] {\bf The Local Index Formula for $SU_q(2)$, ($q=0$)}
\item[5.] {\bf The $\eta$-Cochain}
\item[6.] {\bf Pseudo-differential calculus and the cosphere bundle on ${\rm SU}_q (2)$
, $q \in \, ] 0,1 [$}
\item[7.]{\bf Dimension Spectrum and Residues for ${\rm SU}_q (2)$,  $q \in 
\, ]0,1[$}
\item[8.] {\bf The local index formula for ${\rm SU}_q (2)$, $q \in \, ] 
0,1[$}
\item[9.] {\bf Quantum group Symmetries and Invariant Cyclic Cohomology }
\end{enumerate}

\bigskip
\section{Introduction}

In noncommutative geometry a geometric space is
described from a spectral point of view, as a triple $(\Ac,\Hc ,D)$ consisting of
a $*$-algebra $\Ac$ represented in a Hilbert space $\Hc$ together with an
unbounded selfadjoint operator $D$, with compact resolvent, which interacts
with the algebra in a bounded fashion. This spectral data embodies both
the metric and the differential structure of the geometric space.

\noindent An essential ingredient of the general theory is the Chern character in K-homology
which together with cyclic cohomology and the spectral sequence relating
it to Hochschild cohomology, were defined in 1981 (cf. \cite{Coober},\cite{Coober2},\cite{Coober3}).
The essence of the theory is to allow for computations of differential geometric
nature in the non-commutative framework.

\noindent While basic examples such as the non-commutative tori were
analysed as early as 1980 (cf. \cite{Co80}), the case of the underlying NC-spaces
to quantum groups has been left aside till recently, mainly because of 
the "drop of dimension" which occurs when the deformation parameter $q$
affects non-classical values $q \neq 1$. Thus for instance the Hochschild
dimension of $SU_q(2)$ drops from the classical value $d=3$ to $d=1$
and these NC-spaces seem at first rather esoteric.

\noindent A very interesting spectral triple for  $SU_q(2)$, $q \neq 1$,
has been proposed in \cite{CP}. Thus the algebra $\Ac$ is the algebra of 
functions on  $SU_q(2)$ and the representation in $\Hc$ is the coregular
representation of $SU_q(2)$. The operator $D$ is very simple, 
and is invariant under the action of the quantum group $SU_q(2)$.
(The Anzats proposed in a remark at the end of \cite{CL} provides the right
formula for $\vert D \vert$ but not for the sign of $D$ as pointed out
in \cite{G}).

\noindent Our purpose in this paper is to show that the general theory
developped by Henri Moscovici and the author (cf.\cite{CM}) applies perfectly to
the above spectral triple.

\noindent The power of the general theory comes from 
 general theorems such as the local computation of the analogue of
Pontrjagin classes: {\it i.e.} of the components of the cyclic
cocycle  which is the Chern character of the K-homology class of $D$ and
which
make sense in general. This result allows, using the 
infinitesimal 
calculus, 
to go from local to global in the general framework of spectral triples $(\Ac 
,\Hc,D)$. The notion of locality which is straightforward for classical
spaces is more elaborate in the non-commutative situation and relies
essentially on the non-commutative integral which is the Dixmier trace
in the simplest case and the analogue of the Wodzicki residue in general.
Its validity requires the discreteness of the dimension spectrum, a subset
of $\Cb$ which is an elaboration of the classical notion of dimension.
At an intuitive level this subset is the set of "dimensions", possibly complex, in which
the NC-space underlying the spectral triple manifests itself non-trivially.
At the technical level it is the set of singularities of functions,
\begin{equation}
\z_b (z) = {\rm Trace} \, (b \vert D \vert^{-z}) \qquad {\rm Re}\, z
> p \ , \ b\in \Bc \, .
                                                        \label{eqzet}
\end{equation}
where $ b\in \Bc $ varies in a suitable algebra canonically associated
to the triple and allowing to develop the pseudo-differential 
calculus.

\noindent Our first result is that in the above case of $SU_q(2)$,
the dimension spectrum is {\rm simple} and equal to $\{ 1,2,3 \} \sbs 
\Cb$. Simplicity of the dimension spectrum means that the singularities of the functions (\ref{eqzet})
are at most simple poles. It then follows from the general results of \cite{CM}
that the equality,
\be
\int \!\!\!\!\!\! - P = {\rm Res}_{z=0} \,  \, {\rm Trace}(P\vert
D\vert^{-z})
\ee
defines a  trace on the algebra generated by
$\Ac$, $[D,\Ac]$ and $\vert D \vert^z$, where
$z\in \Cb$.

\noindent Our second result is the explicit computation of this
functional in the above case of $SU_q(2)$. In doing so we shall
also determine the analogue of the cosphere bundle in that 
example and find an interesting space 
$S^*_q$. 
This space is endowed with a one parameter group $\gamma_t$
of automorphisms playing the role of the geodesic flow, and
is intimately related to the product $D^2_{q+} \times D^2_{q-}$, of two NC-two-disks,
while the coproduct gives its relation to $SU_q(2)$.
The formulas computing the residue will be "local"
and very simple, locality by stripping all the expressions from
irrelevant details makes them computable.

\noindent Our third result which is really the main point of the 
paper, is the explicit formula for the local index cocycle, which
owing to the metric dimension $3$ is a priori given by the following cocycle,
\begin{eqnarray}
\label{eq1.70}
\vp_1 (a^0 , a^1) = \int\!\!\!\!\!\!- \ a^0 [D ,a^1 ]\vert D 
\vert^{-1} - \frac{1}{4} \int\!\!\!\!\!\!- \ a^0 \nabla ([D ,a^1 ]) 
\, \vert D \vert^{-3} \\
+ \frac{1}{8} \int\!\!\!\!\!\!- \ a^0 \nabla^2 
([D ,a^1 ]) \, \vert D \vert^{-5} \nonumber 
\end{eqnarray}
and,
\be 
\label{eq1.71}
\vp_3 (a^0 , a^1 , a^2 , a^3) = \frac{1}{12} 
\int\!\!\!\!\!\!- \ a^0 [D ,a^1 ][D , a^2] [D , a^3 ]\, \vert 
D \vert^{-3} \, ,
\ee
where $\nb (T) = [D^2 ,T] \quad  \fl T \
\hbox{operator in} \ \Hc $.
We shall begin by working out the degenerate case $q=0$ with a luxury of details, mainly
to show that the numerical coefficients involved in the above formula
are in fact unique in order to get a (non-trivial) cocycle. The coboundary
involved in the formula (theorem 3) will then be conceptually explained
(in the section "$\eta$-Cochain") and the specific values $\Zc(0)=-\frac{1}{2}$
and $\Zc(-1)=-\frac{1}{12}$ of the Riemann Zeta function
will account for the numerical coefficients encountered in the coboundary.

\noindent We shall then move on to the general case $q \in ]0,1[$ 
 and construct the pseudo-differential calculus on 
${\rm SU}_q (2)$ following the general theory of \cite{CM}. We shall 
determine the algebra of complete symbols by computing the quotient
  by 
smoothing operators. This will give the cosphere bundle $S_q^*$ 
of ${\rm SU}_q (2)$ already mentionned above. The analogue of the 
geodesic flow will give a one-parameter group of
automorphisms $\g_t$ of $C^{\ify} (S_q^*)$. We shall also
construct the restriction morphism $r$ to the 
product of two non-commutative 2-disks,
 
\be
\label{eq4.1}
r : C^{\ify} (S_q^*) \ra C^{\ify} (D_{q_+}^2 \ts D_{q_-}^2)
\ee
We shall then show that the dimension spectrum of ${\rm SU}_q (2)$ 
in the above spectral sense is $\{ 1,2,3 \}$ and compute the 
residues in terms of the symbol $\rho (b) \in C^{\ify} (S_q^*)$
of the operator $b$ of order $0$.
If one lets $\rho (b)^0$ be the component of degree $0$ for the 
geodesic flow $\g_t$, the formulas for the residues are,
$$
\int\!\!\!\!\!\!- \, b \, \vert D \vert^{-3} = (\tau_1 \ot \tau_1)(r\rho 
(b)^0)
$$
$$
\int\!\!\!\!\!\!- \, b \, \vert D \vert^{-2} = (\tau_1 \ot \tau_0 + \tau_0 
\ot \tau_1)(r\rho (b)^0)
$$
$$
\int\!\!\!\!\!\!- \, b \, \vert D \vert^{-1} = (\tau_0 \ot \tau_0)(r\rho 
(b)^0)\, ,
$$
where $r$ is the above restriction map to $D_{q_+}^2 \ts D_{q_-}^2$.
The algebras $C^{\ify} (D_{q \pm}^2)$ are Toeplitz algebras and
as such are extensions of the form,
\be
\label{eq4.2}
0 \longra \Sc \longra C^{\ify} (D^2_{q \pm}) \build\longra_{}^{\s} C^{\ify} (S^1) 
\longra 0
\ee
where the ideal $\Sc$ is isomorphic to the algebra of matrices of rapid 
decay.
The functional $\tau_1$ is the trace obtained by integrating
$\s(x)$ on $S^1$, while  $\tau_0$ 
is a regularized form of the trace on the ideal $\Sc$.
Due to the need of regularization, $\tau_0$ is not a trace
but its Hochschild coboundary (which measures the failure
of the trace property) is easily computed in terms of the 
canonical morphism $\s$.

\noindent A similar long exact sequence, and pair of functionals
$\tau_j$ make sense for $\Ac =C^{\ify} (SU_q(2))$. They are 
invariant under the one parameter group of automorphisms
generated by the derivation $\partial$, which rotates
the canonical generators in opposite ways. Using
this derivation together with the 
second derivative
of $\s(x)$ to define the differential
we then show how to construct a one dimensional cycle
(in the sense of (\cite{Coober2})) whose character
is extremely simple to compute. This shows how to bypass
the shortage of traces on  $\Ac =C^{\ify} (SU_q(2))$
to obtain a significant calculus.

\noindent Our main result (theorem 5) is that the local
formula for the Chern character of the above spectral
triple gives exactly the above cycle, thus completing the
original computation. 
Another quite remarkable point is that the
 cochain whose coboundary
is the difference between the original 
Chern character and the local one is given by the remainders in 
the rational approximation of the logarithmic derivative of the 
Dedekind eta function. The computation of this non-local cochain is very involved.

\noindent One fundamental property of the above spectral triple
is its equivariance (\cite{CP}) under the action of the quantum group $SU_q(2)$.
In the last section we shall use this example to obtain
and explain in general a new
concept
 of quantum group invariance
in cyclic cohomology. 

\bigskip
\section{Operator theoretic Local Index Formula}

\noindent Let   $(\Ac,\Hc,D)$ be a spectral triple.
  The  Fredholm index of the operator  $D$
determines (in the odd case) an  additive map  $K_1 (\Ac) \longra 
 \Zb$ given by the equality
\be
\vp ([u]) = {\rm Index} \, (PuP)  ,  u \in GL_1 (\Ac)                                                 \label{eq:(3.3)}
\ee
where $P$ is the projector  $P = \frac{1+F}{2}$, $F = {\rm Sign} \, (D)$.

\smallskip
 
 \noindent This map is computed by the pairing
of  $K_1 (\Ac)$ with the following  cyclic cocycle
\be
\tau (a^0 ,\ldots ,a^n) ={\rm Trace} \, (a^0 [F,a^1] \ldots
[F,a^n]) \qquad \fl \, a^j \in \Ac
                                                \label{eq:(3.4)}
\ee
where $F= \hbox{Sign} \ D$ and  we assume that the dimension $p$ of our space is 
finite,
which means that the characteristic values $\mu_k$ of $(D+i)^{-1}$ 
decay like $k^{-1/p}$, also $n\geq p$ is an odd integer. 
There are similar formulas involving the grading $\gamma$ in the even case.

\noindent The cyclic cohomology $HC^n (\Ac)$ is defined as the cohomology
of the complex of cyclic cochains, i.e. those satisfying

\be
\label{cyclic}
\psi (a^1 , \ldots ,a^n ,a^0 ) = (-1)^n \ \psi (a^0 ,\ldots , a^n)  \ , \qquad \fl
a^j \in \Ac \ , 
\ee
under the coboundary operation $b$ given by:
\begin{eqnarray}
&& (b\psi)(a^0 ,\ldots ,a^{n+1}) =
\sum_0^n (-1)^j \ \psi (a^0 ,\ldots ,a^j \ a^{j+1} ,\ldots ,a^{n+1})\\
&&+ (-1)^{n+1}
\ \psi (a^{n+1} \ a^0 ,\ldots ,a^n ) \ , \quad \fl a^j \in \Ac \ .
\nonumber 
\end{eqnarray}

\noindent Equivalently, $HC^n (\Ac)$ can be described in terms of the second filtration of
 the $(b,B)$ bicomplex of arbitrary (non cyclic) cochains on $\Ac$, where 
$B:C^m \rightarrow C^{m-1}$ is given by
\begin{eqnarray}
&&(B_0  \vp) (a^0 ,\ldots ,a^{m-1}) = \vp (1,a^0 ,\ldots ,a^{m-1}) - (-1)^m 
\vp (a^0 ,\ldots ,a^{m-1},1) \nonumber \\
&& B=AB_0 , \quad (A\psi) (a^0 ,\ldots ,a^{m-1}) = \sum (-1)^{(m-1)j} \ \psi
(a^j ,\ldots ,a^{j-1})  
\end{eqnarray}

\noindent To an $n$-dimensional cyclic cocycle $\psi$ one associates the
$(b,B)$ cocycle $\vp \in Z^p (F^q \ C)$, $n=p-2q$ given by
\be
\label{cyctobb}
(-1)^{[n/2]} \ (n!)^{-1} \ \psi = \vp_{p,q} 
\ee
where $\vp_{p,q}$ is the only non zero component of $\vp$.

\noindent Given a
spectral triple $(\Ac ,\Hc ,D)$, with $D^{-1} \in \Lc^{(p,\ify)}$, the
precise normalization for its
Chern character in cyclic cohomology is obtained from the following cyclic
cocycle $\tau_n$, $n\geq p$, $n$ odd,
\be
\label{chern}
\tau_n (a^0 ,\ldots ,a^n) = \lb_n \ {\rm Tr}' \left( a^0 [F,a^1] \ldots [F,a^n]
\right)\ , \qquad \fl a^j \in \Ac \ , 
\ee
where $F={\rm Sign} D$, $\lb_n = \sqrt{2i} \ (-1)^{{n(n-1) \over 2}} \ \G
\left( {n\over 2} +1 \right)$ and
\be
{\rm Tr}' (T) = {1\over 2} \ {\rm Trace} (F(FT+TF)) \, 
\ee 
If one wants to 
regard the cocycle $\tau_n$ of
(\ref{chern}) as a cochain of the $(b,B)$ bicomplex, one takes (\ref{cyctobb})
into account and use instead of $\lb_n$, the normalization constant
$\mu_n = (-1)^{[n/2]} \ (n!)^{-1} \ \lb_n = \sqrt{2i} \ {\G \left( {n\over 2}
+1\right) \over n!}$.

\smallskip

\noindent It is difficult to compute the cocycle $\tau_n$ in general because the formula
(\ref{chern}) involves the ordinary trace instead of the local trace
${\int \!\!\!\!\! -}$ and it is crucial to obtain a local form of the above 
cocycle.

 \noindent In \cite{book} we obtained the following general formula for 
the Hochschild cohomology class of
$\tau_n$ in terms of the Dixmier trace :
\be
\vp_n (a^0 ,\ldots ,a^n) = \lb_n \ {\rm Tr}_{\om} \left( a^0 [D,a^1] \ldots
[D,a^n] \ \vert D \vert^{-n} \right) \ , \quad \fl a^j \in \Ac \ . 
\ee
The problem of finding a local formula for the 
{\it cyclic cohomology} Chern character, i.e. for the 
class of $\tau_n$ is solved by a general formula \cite{CM} which is expressed in terms of the $(b,B)$
bicomplex and which we now explain.

\noindent Let us make the following regularity hypothesis
on  $(\Ac ,\Hc ,D)$
\begin{equation}
a \ \hbox{and } \ [D,a] \ \in \ \cap  \, {\rm Dom} \, \d^k , \ \fl \, a\in
\Ac                                                     \label{eq:(3.1)}
\end{equation}
where $\d$ is the derivation $\d(T) = [\vert D \vert,T]$ for any operator $T$.

\smallskip

\noindent We let $\Bc$ denote the algebra generated by  $\d^k (a)$,
$\d^k ([D,a])$. 
The usual notion of {\it dimension}  of a space is replaced by the
{\it dimension spectrum } which is a subset of $\Cb$.
\noindent The precise definition of the dimension spectrum is the subset
$\Si \sbs \Cb$ of singularities of the analytic functions
\begin{equation}
\z_b (z) = {\rm Trace} \, (b \vert D \vert^{-z}) \qquad {\rm Re} z
> p \ , \ b\in \Bc \, .
                                                        \label{eq:(3.2)}
\end{equation}
Note that $D$ may have a non-trivial kernel 
 so that $\vert D \vert^{-s}$
is ill defined there. However the kernel of $D$ is finite 
dimensional and the poles and residues of the above function
are independent
of the arbitrary choice of a non-zero positive value $\vert D \vert=\ve$
on this kernel.
The dimension spectrum of a manifold $M$ consists of relative integers less than $n=\dim M$; it is simple.  Multiplicities appear for
singular manifolds. Cantor sets provide examples of  complex points $z
\notin \Rb$ in the dimension  spectrum.

\noindent We assume  that $\Si$ is discrete and simple, i.e. that
$\z_b$ can be extended to  $\Cb / \Si$  with simple poles in  $\Si$.
In fact the hypothesis only matters in a neighborhood of $\{z,{\rm Re}(z)\geq 
0 \}$.

\smallskip

\noindent  Let $(\Ac ,\Hc ,D)$ be a spectral triple satisfying 
the hypothesis
(\ref{eq:(3.1)}) and (\ref{eq:(3.2)}).
\smallskip

\noindent  We shall use the following notations:
$$
\nb (a) = [D^2 ,a] \quad ; \quad a^{(k)} = \nb^k (a) \ , \quad \fl a \
\hbox{operator in} \ \Hc \ . 
$$

\noindent The local index theorem is the following,  
\cite{CM}:

\bigskip

\noindent {\bf Theorem 1.}  \begin{enumerate}
\item {\it The equality
\[
{\int \!\!\!\!\!\! -} P = {\rm Res}_{z=0} \,  \, {\rm Trace}(P\vert
D\vert^{-z})
\]
defines a  trace on the algebra generated by
$\Ac$, $[D,\Ac]$ and $\vert D \vert^z$, where}
$z\in \Cb$.

\item
{\it There is  only  a finite number of non--zero terms in the following formula which
defines the odd components $(\vp_n)_{n=1,3,\ldots}$
of a cocycle in the bicomplex $(b,B)$ of $\Ac$,
 \[
\vp_n (a^0 ,\ldots ,a^n) = \sum_k c_{n,k}
{\int \!\!\!\!\!\! -} a^0 [D,a^1]^{(k_1)} \ldots
[D,a^n]^{(k_n)} \, \vert D \vert^{-n-2\vert k\vert}
\qquad \fl \, a^j \in \Ac
 \]
where $k$ is a multi-index}, $\vert k \vert = k_1 +\ldots + k_n$,
\[
c_{n,k} =
(-1)^{\vert k \vert} \, \sqrt{2i} (k_1 ! \ldots k_n !)^{-1} \, ((k_1 +1) \ldots (k_1 + k_2 + \ldots + k_n
+n))^{-1} \, \G \left( \vert k \vert + {n\over 2}
\right).
\]

\item {\it The pairing of the cyclic cohomology class $(\vp_n) \in HC^* (\Ac)$
with $K_1 (\Ac)$ gives the Fredholm index of $D$ with coefficients in}
$K_1(\Ac)$.
\end{enumerate}

\bigskip

\noindent For the normalization of the pairing
between $HC^*$ and  $K(\Ac)$ see \cite{book}.
In the even case, i.e. when $\Hc$ is $\Zb /2$ graded by $\g$,
 \[
\g = \g^*, \ \ \g^2 =1, \ \ \g a = a \g \quad \fl \, a \in \Ac, \ \g
D = -D\g,
\]
there is an analogous formula for a cocycle
$(\vp_n)$, $n$ even, which gives the Fredholm index of $D$
with  coefficients in  $K_0$. However,
$\vp_0$ is not expressed  in terms of the residue
${\int \!\!\!\!\! -}$ because the character can be non-trivial
for a finite dimensional  $\Hc$, in which case all
residues vanish.

\noindent To give some concreteness to this general result we 
shall undertake the computation in an example, that of the 
quantum group ${\rm SU}_q (2)$. Its original interest is that 
it lies rather far from ordinary manifolds and is thus a good test
case for the general theory.

\section{Dimension Spectrum of ${\rm SU}_q (2)$: Case $q=0$.}

Let $q$ be a real number $0 \leq q<1$.
We start with the presentation of the algebra of coordinates on the quantum group
${\rm SU}_q (2)$ in the form, 

\be
\label{eq1}
\a^* \a + \b^* \b = 1 \, , \ \a \a^* + q^2 \b \b^* = 1 \, , \ 
\a \b = q \b \a \, , \ \a \b^* = q \b^* \a \, , \ \b \b^* = 
\b^* \b \, .
\ee \noindent Let us recall the notations for the standard representation 
of that algebra. One lets $\Hc$ be the Hilbert space with orthonormal basis 
$e_{ij}^{(n)}$ where $n \in \frac{1}{2} \, \Nb$ varies among half-integers while 
$i,j \in \{ -n , -n+1 , \ldots , n \}$.

\noindent Thus the first elements are,
$$
e_{00}^{(0)} \, , \ e_{ij}^{(1/2)} \, , \ i,j \in \left\{ -  \frac{1}{2} , \frac{1}{2} \right\} , \ldots
$$
The following formulas define a unitary 
representation in $\Hc$,
\be
\label{eq2}
\a e_{ij}^{(n)} = a_+ (n,i,j) \, e_{i-\frac{1}{2} , j - 
\frac{1}{2}}^{\left( n + \frac{1}{2} \right)} + a_- (n,i,j) 
\, e_{i-\frac{1}{2} , j - \frac{1}{2}}^{\left( n - 
\frac{1}{2} \right)}
\ee
$$
\b e_{ij}^{(n)} = b_+ (n,i,j) \, e_{i+\frac{1}{2} , j - 
\frac{1}{2}}^{\left( n + \frac{1}{2} \right)} + b_- (n,i,j) 
\, e_{i+\frac{1}{2} , j - \frac{1}{2}}^{\left( n - 
\frac{1}{2} \right)}
$$
where the explicit form of $a_{\pm}$ and $b_{\pm}$ 
 is,
\be
\label{eq25}
a_+ (n,i,j) = q^{2n+i+j+1} \, \frac{(1-q^{2n-2j+2})^{1/2} 
(1-q^{2n-2i+2})^{1/2}}{ (1-q^{4n+2})^{1/2} (1-q^{4n+4})^{1/2}}
\ee

$$
a_- (n,i,j) = \frac{(1-q^{2n+2j})^{1/2} 
(1-q^{2n+2i})^{1/2}}{(1-q^{4n})^{1/2} (1-q^{4n+2})^{1/2}}
$$
and
\be
\label{eq26}
b_+ (n,i,j) = -q^{n+j} \, \frac{(1-q^{2n-2j+2})^{1/2} 
(1-q^{2n+2i+2})^{1/2}}{(1-q^{4n+2})^{1/2} (1-q^{4n+4})^{1/2}}
\ee

$$
b_- (n,i,j) = q^{n+i} \, \frac{(1-q^{2n+2j})^{1/2} 
(1-q^{2n-2i})^{1/2}}{(1-q^{4n})^{1/2} (1-q^{4n+2})^{1/2}} \, .
$$
Note that $a_-$ does vanish if $i=-n$ or $j=-n$, which gives 
meaning to $a_- (n,i,j) \, e_{i-\frac{1}{2} , 
j-\frac{1}{2}}^{\left( n-\frac{1}{2}\right)}$ for these values 
while $i-\frac{1}{2} \notin \left[ - \left( n-\frac{1}{2} 
\right) , n - \frac{1}{2} \right]$ or $j - \frac{1}{2} \notin \left[ - \left( n-\frac{1}{2} 
\right) , n - \frac{1}{2} \right]$. Similarly $b_-$ vanishes for $j=-n$ or $i=n$.
\bigskip

\noindent  Let now as in (\cite{CP}), $D$ be the diagonal operator in $\Hc$ given by,
\be
\label{eq2.18}
D (e_{ij}^{(n)}) = (2 \, \d_0 (n-i)-1) \, 2n \,\, e_{ij}^{(n)}
\ee
where $\d_0 (k) = 0$ if $k \ne 0$ and $\d_0 (0) = 1$. It follows from 
\cite{CP} that the triple
\be
\label{eq2.19}
(\Ac , \Hc , D)
\ee
is a spectral triple.

\smallskip
\noindent  In order to simplify we start the discussion  with the case $q=0$.
We then have the simpler formulas,
\begin{eqnarray}
\label{eq1.1}
&&a_+ (n,i,j) = 0 \\
&&a_- (n,i,j) = 0 \ {\rm if} \ i=-n \ {\rm or} \ j = -n 
\nonumber \\
&&a_- (n,i,j) = 1 \ {\rm if} \ i \ne -n \ {\rm and} \ j \ne -n 
\nonumber
\end{eqnarray}

\begin{eqnarray}
\label{eq1.2}
&&b_+ (n,i,j) = 0 \ {\rm if} \ j \ne -n \\
&&b_+ (n,i,j) = -1 \ {\rm if} \ j = -n \nonumber \\
&&b_- (n,i,j) = 0 \ {\rm if} \ i\ne -n \ {\rm or} \ j = -n 
\nonumber \\
&&b_- (n,i,j) = 1 \ {\rm if} \ i=-n , j \ne -n \, . \nonumber
\end{eqnarray}

\smallskip 
 \noindent  Thus for $q=0$ the operators $\a$ and $\b$ in $\Hc$ are given by,
\be
\label{eq2.1}
\a e_{ij}^{(n)} = e_{i - \frac{1}{2} , j - \frac{1}{2}}^{\left( n - 
\frac{1}{2} \right)} \qquad \hbox{if} \quad i > -n , j > -n
\ee
and $\a e_{ij}^{(n)} = 0$ if $i = -n$ or $j=-n$.
\be
\label{eq2.2}
\b e_{ij}^{(n)} = 0 \qquad \hbox{if} \quad i \ne -n \quad \hbox{and} 
\quad j \ne -n
\ee
\be
\label{eq21.3}
\b e_{-n,j}^{(n)} = e_{-\left( n - \frac{1}{2} \right)  , j - 
\frac{1}{2}}^{\left( n - \frac{1}{2} \right)} \qquad \hbox{if} \quad j 
\ne -n
\ee
and
\be
\label{eq21.4}
\b e_{i,-n}^{(n)} = - e_{i+ \frac{1}{2} , - \left( n+ \frac{1}{2} 
\right)}^{\left( n + \frac{1}{2} \right)} \, .
\ee
By construction $\b \b^* = \b^* \b$ is the projection $e$ on the subset 
$\{ i = -n$ or $j = -n\}$ of the basis.

\smallskip

\noindent Also $\a$ is a partial isometry with initial support $1-e$ and 
final support $1 = \a \a^*$. The basic relations between $\a$ and $\b$ 
are,
\be
\label{eq2.5}
\a^* \a + \b^* \b = 1 \, , \ \a \a^* = 1 \, , \ \a \b = \a \b^* = 0 \, , 
\ \b \b^* = \b^* \b \, .
\ee
For $f \in C^{\ify} (S^1)$, $f = \sum \wh f_n \, e^{in\t}$, we let
\be
\label{eq2.6}
f(\b) = \sum_{n > 0} \wh f_n \, \b^n + \sum_{n < 0} \wh f_n \, 
\b^{*(-n)} + \wh f_0 \, e
\ee
and the map $f \ra f(\b)$ gives a (degenerate) representation of 
$C^{\ify} (S^1)$ in $\Hc$.

\smallskip

\noindent Now let $\Ac$ be the linear space of sums,
\be
\label{eq2.7}
a = \sum_{k,\ell \geq 0} \a^{*k} f_{k\ell} (\b) \, \a^{\ell} + 
\sum_{\ell \geq 0} \lb_{\ell} \, \a^{\ell} + \sum_{k > 0} \lb'_k \, 
\a^{*k}
\ee
where $\lb$ and $\lb'$ are sequences (of complex numbers) of rapid decay 
and $(f_{k\ell})$ is a sequence of rapid decay with values in $C^{\ify} 
(S^1)$.

\smallskip

\noindent We let $A$ be the $C^*$ algebra in $\Hc$ generated by $\a$ and 
$\b$. 
\medskip

\noindent {\bf Proposition 1.} {\it The subspace $\Ac \sbs A$ is a 
subalgebra stable under holomorphic functional calculus.}

\medskip

\noindent {\bf Proof.} Let $\s$ be the linear map from $\Ac$ to 
$C^{\ify} (S^1)$ given by,
\be
\label{eq2.8}
\s (a) = \sum_{\ell \geq 0} \lb_{\ell} \, u^{\ell} + \sum_{k > 0} \lb'_k 
\, u^{-k}
\ee
where $u = e^{i\t}$ is the generator of $C^{\ify}  (S^1)$. Let $\Jc = 
{\rm Ker} \, \s$. For $a \in \Jc$ one has $a = \sum \a^{*k} f_{k\ell} 
\, \a^{\ell}$ and the equality,
\be
\label{eq2.9}
\a^{*k} f_{k\ell} \, \a^{\ell} \a^{*k'} g_{k'\ell'} \, \a^{\ell'} = 
\d_{\ell,k'} \, \a^{*k} f_{k\ell} \, g_{k'\ell'} \, \a^{\ell'}
\ee
shows that $\Jc$ is an algebra and is isomorphic to the topological 
tensor product
\be
\label{eq2.10}
C^{\ify} (S^1) \ot \Sc = C^{\ify} (S^1 , \Sc)
\ee
where $\Sc$ is the algebra of matrices of rapid decay.

\smallskip

\noindent Since $\Sc$ is stable under holomorphic functional calculus 
(h.f.c.) in its norm closure $\Kc$ (the $C^*$ algebra of compact 
operators), it follows from (\ref{eq2.10}) that $\Jc$ is stable under 
h.f.c. in its norm closure $\ov\Jc \sbs A$.

\smallskip

\noindent The equalities $\a f(\b) = 0$ $\fl f \in C^{\ify} (S^1)$ and 
$\a \a^* = 1$ show that $\Jc$ is stable under left multiplication by 
$\a^*$ and $\a$. It follows using (\ref{eq2.5}) that $\Ac$ is an 
algebra, $\Jc$ a two sided ideal of $\Ac$ and that one has the exact 
sequence,
\be
\label{eq2.11}
0 \longra \Jc \longra \Ac \build\longra_{}^{\s} C^{\ify} (S^1) \longra 0 
\, .
\ee
By construction $\Ac$ is dense in $A$. Let us check that it is stable 
under h.f.c. in $A$. Let $a \in \Ac$ be such that $a^{-1} \in A$. Let us
 show that $a^{-1} \in \Ac$.

\smallskip

\noindent Let $\partial_{\a}$ be the derivation of $\Ac$ given by,
\be
\label{eq2.12}
\partial_{\a} \a = \a \, , \ \partial_{\a} \b = 0 \, .
\ee
The one parameter group $\exp (i t \partial_{\a})$ of automorphisms of $\Ac$ 
is implemented by unitary operators in $\Hc$ (cf.(\ref{eq1.43}) below)
and extends to $A$. 
Moreover $\Ac$ is dense in the domain,
\be
\label{eq2.13}
{\rm Dom} \, \partial_{\a}^j = \{ x \in A \, ; \ \partial_{\a}^j x \in A \}
\ee
in the graph norm.

\smallskip

\noindent Since $a^{-1} \in {\rm Dom} \, \partial_{\a}^j$ 
we can, given any $\ve > 0$, find $b \in \Ac$ such that,
\be
\label{eq2.14}
\Vert \partial_{\a}^j (b-a^{-1}) \Vert < \ve \qquad j = 0,1,2 \, .
\ee

\smallskip

\noindent Thus, given $\ve > 0$, we can find $b \in \Ac$ such that, with 
$x=ab$,
\be
\label{eq2.15}
\Vert \partial_{\a}^j (x-1) \Vert < \ve \qquad j = 0,1,2 \, .
\ee
For $\ve$ small enough it follows that if we let $\s (x^{-1})_n^{\wdg} $ be the  Fourier coefficients of $\s (x^{-1})$, 
\be
\label{eq2.16}
c = \sum_{n \geq 0} \s (x^{-1})_n^{\wdg} \, \a^n + \sum_{n < 0} \s 
(x^{-1})_n^{\wdg} \, \a^{*-n}
\ee
is an element of $\Ac$,  invertible in  $A $, such that,
\be
\label{eq2.17}
\s (c) = \s (x^{-1}) \, .
\ee
(Since one controls $n^2 \s (x^{-1})_n^{\wdg}$ from $\Vert \partial_{\a}^j 
(x^{-1} - 1) \Vert$.)

\smallskip

\noindent Thus $\s (xc) = 1$ and since $xc$ is invertible in $A$ (by 
(\ref{eq2.15}), (\ref{eq2.17})) and $1-xc \in \Jc$ the stability of 
$\Jc$ under h.f.c. shows that $(xc)^{-1} = y \in \Ac$. Then $abcy = 1$  and $a^{-1} = bcy \in \Ac$. $\Box$

\smallskip

\noindent Our next result determines the dimension spectrum of the 
spectral triple $(\Ac , \Hc , D)$ defined above in 
(\ref{eq2.19}),

\bigskip

\noindent {\bf Theorem 2.} {\it The dimension spectrum of the spectral 
triple $(\Ac , \Hc , D)$ is {\rm simple} and equal to $\{ 1,2,3 \} \sbs 
\Cb$.}

\bigskip

\noindent Thus we let $\Bc$ be the algebra generated by the
\be
\label{eq2.20}
\d^k (a) \, , \ \d^k ([D,a]) \, , \ a \in \Ac \, , \ k \in \Nb
\ee
where $\d$ is the unbounded derivation of $\Lc (\Hc)$ given by the 
commutator with $\vert D \vert$,
\be
\label{eq2.21}
\d (T) = \vert D \vert \, T - T \, \vert D \vert \, .
\ee
(It is part of the statement that the elements in (\ref{eq2.20}) are in 
the domain of $\d^k$.)

\smallskip

\noindent For $b \in \Bc$ we consider the function,
\be
\label{eq2.22}
\Zc_b (s) = {\rm Trace} \, (b \, \vert D \vert^{-s})
\ee
where we take care of the eigenvalue $D=0$ by replacing $\vert D \vert$ 
by an arbitrary $\ve > 0$ there. The statement of the theorem is that 
all the functions $\Zc_b (s)$ which are a priori only defined for ${\rm 
Re} \, (s) > 3$, do extend to meromorphic function on $\Cb$ and only 
admit {\it simple} poles at the 3 points $\{1,2,3\} \sbs \Cb$.

\smallskip

\noindent To prove it we shall first describe the algebra $\Bc$. We let,
\be
\label{eq2.23}
F = {\rm Sign} \, D
\ee
so that $F = 2P-1$ where $P$ is the orthogonal projection on the subset 
$\{ i = n \}$ of the basis. Concerning the generator $\a$ one has,
\be
\label{eq2.24}
\d (\a) = - \a \, , \ \d (\a^*) = \a^* \, , \ [F,\a] = 0 \, .
\ee
It follows that $[D,\a^*] = F \d (\a^*) = F \a^* = \a^* F$. Thus $\a \, 
[D,\a^*] = \a \a^* F = F$,
\be
\label{eq2.25}
F = \a \, [D , \a^*] \, .
\ee
This shows that $F \in \Bc$.

\smallskip

\noindent Concerning the generator $\b$ one has
\be
\label{eq2.26}
\d (\b) = \b K \, , \ \d (\b^*) = - K \b^*
\ee
where $K$ is the multiplication operator,
\be \label{eq2.27}
K (e_{ij}^{(n)}) = k(n,i,j) \, e_{ij}^{(n)}
\ee
with
\be
\label{eq2.28}
k (n,i,j) = 0 \qquad\hbox{unless} \ i = -n \quad \hbox{or} \quad j = 
-n
\ee
$$
k(n,-n,j) = -1 \qquad \hbox{if}   \ j \ne -n
$$
$$
k(n,i,-n) = 1 \, . 
$$
Thus the support of $K$ is $e = \b^* \b$ and,
\be
\label{eq2.29}
K^2 = e \, .
\ee
We let $e_0 = \frac{1}{2} \, (K-\b K \b^*)$. It is the orthogonal 
projection on the subset of the basis $\{ i=-n$ and $j=-n\}$. For each 
$m$ one lets,
\be \label{eq2.30}
e_m = \b^m e_0 \, \b^{*m}
\ee and the $e_m$ are pairwise orthogonal projections such that,
\be
\label{eq2.31}
\sum_{m \in \Zb} e_m = e \, . \ee
We let $\Lc$ be the algebra of double sums with rapid decay,
\be
\label{eq2.32}
\Lc = \left\{ \sum \lb_{n,m} \, \b^n e_m \, ; \ \lb \in \Sc \right\}
\ee
(where $\b^{-\ell} = \b^{*\ell}$ for $\ell > 0$).

\smallskip

\noindent One has $\d (e_n) = 0$, $[K,\b ]=2 e_0 \b$, $K e_m= $ sign$(m) e_m$
  and using (\ref{eq2.26}),
\be
\label{eq2.33}
\d (\b^n) = n \b^n K \qquad \hbox{modulo} \ \Lc \, .
\ee
Thus $\Lc$ is invariant (globally) under $\d$. Also for any $f \in 
C^{\ify} (S^1)$ one has
\be
\label{eq2.34}
[K , f(\b)] \in \Lc
\ee
and the algebra $B_0$,
\be
\label{eq2.35}
B_0 = \{ f_0 (\b) + f_1 (\b) \, K + h \, ; \ f_j \in C^{\ify} (S^1) \, , 
\ h \in \Lc \}
\ee
is stable by the derivation $\d$.

\smallskip

\noindent A similar result holds if we further adjoin the operator
\be
\label{eq2.36}
F_1 = eF = Fe \, .
\ee
Indeed $e'_0 = \frac{1}{2} \, (F_1 - \b F_1 \b^*)$ is the projection on 
the element  $\{ e_{0,0}^{(0)} \}$  of the basis and the $e'_n = \b^n e'_0 
\b^{*n}$ are pairwise orthogonal projections on the one dimensional 
subspaces spanned for $n > 0$ by $e_{n/2 , -n/2}^{(n/2)}$ and for $n < 
0$, $n=-k$, by $e_{-k/2 , k/2}^{(k/2)}$.

\smallskip

\noindent We let,
\be
\label{eq2.37}
\Lc' = \left\{ \sum \lb_{n,m} \, \b^n e'_m \, ; \ \lb \in \Sc \right\} \, 
.
\ee
One has, \be \label{eq2.38}
[F,f(\b)] \in \Lc' \qquad \fl f \in C^{\ify} (S^1)
\ee
and $F_1 \Lc' = \Lc' F_1 = \Lc'$.

\smallskip

\noindent Also $e'_n \leq e_n$ for each $n$ so that $\Lc \Lc' \sbs \Lc'$ 
and $\Lc' \Lc \sbs \Lc'$ which shows that the sum,
\be
\label{eq2.39}
\Lc'' = \Lc + \Lc'
\ee
is an algebra.

\smallskip

\noindent Thus the algebra generated in $e \Hc$ by the $\d^k (f)$, $f \in 
C^{\ify} (S^1)$ and $F_1$, is contained in the algebra $B_1$
\be
\label{eq2.40}
B_1 = \{ f_0 + f_1 K + f_2 F_1 + h \, ; \ f_j \in C^{\ify} (S^1) \, , \ h 
\in \Lc'' \} \, .
\ee
Note that $F_1 K = 1+2(K-F_1)$ so that we do not need terms in $F_1 K$.

\smallskip

\noindent We then let $B$ be the algebra of double sums,
\be
\label{eq2.41}
B = \left\{ \sum \a^{*k} b_{k\ell} \, \a^{\ell} + A_0 + A_1 F \right\}
\ee
where $b_{k\ell} \in B_1$ and the sequence $(b_{k\ell})$ is of rapid 
decay while $A_0 , A_1$ are sums of rapid decay of the form,
$$
A = \sum_{\ell \geq 0} a_{\ell} \, \a^{\ell} + \sum_{k > 0} a_{-k} \, 
\a^{*k} \, .
$$
Since $F$ commutes with $\a$ and $\a^*$ it commutes with $A_j$. Thus one 
checks that $B$ is an algebra, that it is stable under $\d$ and contains 
both $F$ and $\Ac$, thus it contains $\Bc$.

\smallskip

\noindent Let then $b \in B$ and consider the function,
\be
\label{eq2.42}
\Zc_b (s) = {\rm Trace} \, (b \, \vert D \vert^{-s})
\ee
which is well defined for ${\rm Re} \, (s) > 3$.

\smallskip
\noindent There is a natural bigrading corresponding to the 
degrees in $\a$ and $\b$. It is implemented by the following action of $\Tb^2$ in $\Hc$, 
\be
\label{eq1.43}
V(u,v) \, e_{k,\ell}^{(n)} = \exp i (-u (k+\ell) + v(k-\ell)) \, 
e_{k,\ell}^{(n)} \, .
\ee
Note that both $k+\ell$ and $k-\ell$ are integers so that one gets an 
action of $\Tb^2$.

\smallskip

\noindent The indices $k,\ell$ are transformed to $k - \frac{1}{2}$, 
$\ell - \frac{1}{2}$ by $\a$, so that
$$
V(u,v) \, \a (e_{k,\ell}^{(n)}) = \exp i (-u (k+\ell -1) + v 
(k-\ell)) \, \a (e_{k,\ell}^{(n)})
$$
and we get,
\be
\label{eq1.44}
V(u,v) \, \a \, V (-u , -v) = e^{iu} \a \, .
\ee
The indices $k,\ell$ are transformed to $k+\frac{1}{2}$, $\ell - 
\frac{1}{2}$ by $\b$ and,
$$
V(u,v) \, \b (e_{k,\ell}^{(n)}) = \exp i (-u(k+\ell) + v(k-\ell+1)) 
\, \b (e_{k,\ell}^{(n)})
$$
so that, \be
\label{eq1.45}
V(u,v) \, \b \,V (-u,-v) = e^{iv} \b \, .
\ee
Moreover, since $V$ is a multiplication operator it commutes with 
$\vert D \vert$, $D$, and $F=2P-1$.

\smallskip

\noindent Using the restriction of this bigrading to $B$ 
(which gives bidegree $(0,0)$ 
for diagonal operators, $(1,0)$ for $\a$ and $(0,1)$ for $\b$) one checks 
that homogeneous elements of bidegree $\ne (0,0)$ satisfy $\Zc_b (s) 
\equiv 0$, thus one can assume that $b$ is of bidegree $(0,0)$.

\smallskip

\noindent Any $b \in B^{(0,0)}$ is of the form,
\be
\label{eq2.43}
b = \sum \a^{*k} \, b_k \, \a^k + a_0 + a_1 F
\ee
where $a_0 , a_1$ are {\it scalars} and $(b_k)$ is a sequence of rapid 
decay with $b_k \in B_1^{(0,0)}$. Elements $c$ of $B_1^{(0,0)}$ are of 
the form,
\be
\label{eq2.44}
c = \lb_0 + \lb_1 K + \lb_2 F_1 + h
\ee where $\lb_j$ are {\it scalars} and $h \in \Lc''^{(0,0)}$. Finally 
elements of $\Lc''^{(0,0)}$ are of the form,
\be
\label{eq2.45}
h = \sum h_n e_n + \sum h'_m e'_m
\ee
where $(h_n)$ and $(h'_m)$ are {\it scalar} sequences of rapid decay.

\smallskip

\noindent The equality,
\be
\label{eqperm}
\vert D \vert^z \, \a^{*k} = \a^{*k} (\vert D \vert + k)^z \qquad z \in 
\Cb
\ee
is checked directly ($k \geq 0$).

\smallskip

\noindent Using $\a^k \a^{*k} = 1$ it follows that with $b$ as in (\ref{eq2.43}),
\be
\label{eq2.47}
{\rm Trace} \, (b \, \vert D \vert^{-s}) = {\rm Trace} \, ((a_0 + a_1 F) 
\, \vert D \vert^{-s}) + \sum_{k \geq 0} {\rm Trace} \, (b_k (\vert D 
\vert + k)^{-s}) \, .
\ee
Now for $h$ as in (\ref{eq2.45}) one has
$$
{\rm Trace} \, (h \, \vert D \vert^{-s}) = \sum h_n {\rm Trace} \, (e_n 
\, \vert D \vert^{-s}) + \sum h'_m {\rm Trace} \, (e'_m \, \vert D  \vert^{-s}) \, . $$
Moreover
$$
{\rm Trace} \, (e_n \, \vert D \vert^{-s}) = \sum_{\ell = 0}^{\ify} 
\frac{1}{(\vert n \vert + \ell)^s} = \Zc (s) - \left( \sum_0^{\vert n 
\vert - 1} \frac{1}{r^s} \right) \, .
$$
But $\sum h_n \left( \build\sum_{0}^{\vert n \vert - 1} \frac{1}{r^s} 
\right) = \rho_1 (s)$ is a holomorphic function of $s \in \Cb$ and 
similarly, since ${\rm Trace} \, (e'_m \, \vert D \vert^{-s}) = 
\frac{1}{\vert m \vert^s}$, the function $\rho_2 (s) = \sum h'_m {\rm 
Trace} \, (e'_m \, \vert D \vert^{-s})$ is holomorphic in $s \in \Cb$. 
Thus modulo holomorphic functions one has,
\be
\label{eq2.48}
{\rm Trace} \, (h \, \vert D \vert^{-s}) \sim \left( \sum h_n \right) \Zc 
(s) \, .
\ee
Next, 
\be
\label{e}
{\rm Trace} \, (e \, \vert D \vert^{-s}) = \sum_0^{\ify} \frac{2n+1}{n^s} 
= 2 \,\Zc (s-1) + \Zc (s)+ \ve^{-s}
\ee
$$
{\rm Trace} \, (K \, \vert D \vert^{-s}) = \sum_{n \in \Zb} \ \sum_{\ell 
\geq 0} \ \frac{{\rm sign} (n)}{(\vert n \vert + \ell)^s} = \Zc (s)+ \ve^{-s}
$$
and with $F = 2P-1$ we also have,
\be
\label{F}
{\rm Trace} \, (e P \, \vert D \vert^{-s}) = \,\Zc (s) \, + \ve^{-s}.
\ee
Thus, with $c$ as in (\ref{eq2.44}) we get,
\be
\label{eq2.49}
\Zc_c (s) = \lb \, \Zc (s-1) + \mu \, \Zc (s) + \rho (s)
\ee
where $\lb , \mu$ are scalars and $\rho$ is a holomorphic function of $s 
\in \Cb$.

\smallskip

\noindent A similar result holds for
$$
\sum_{k \geq 0} {\rm Trace} \, (b_k (\vert D \vert + k)^{-s}) \, .
$$
For instance one rewrites the double sum
$$
\sum h_{n,k} {\rm Trace} \, (e_m (\vert D \vert + k)^{-s}) = 
\sum_{n,k,\ell} h_{n,k} \frac{1}{(\vert n \vert + k + \ell)^s}
$$
as
$$
\sum_m \left( \sum_{\vert n \vert +k \leq m} h_{n,k} \right) 
\frac{1}{m^s} = a \, \Zc(s) + \rho (s)
$$
where $a = \sum h_{n,k}$ and $\rho$ is holomorphic in $s \in \Cb$. 

\smallskip

\noindent Finally
$$
{\rm Trace} \, (P \, \vert D \vert^{-s}) = \sum_0^{\ify} 
\frac{(n+1)}{n^s} = \Zc (s-1) + \Zc (s)
$$
and
$$
{\rm Trace} \, (\vert D \vert^{-s}) = \sum_0^{\ify} \frac{(n+1)^2}{n^s} = 
\Zc (s-2) + 2 \, \Zc (s-1) + \Zc (s) \, .
$$
Thus we conclude that for any $b \in B$ one has
\be
\label{eq2.50}
\Zc_b (s) = \lb_3 \, \Zc (s-2) + \lb_2 \, \Zc (s-1) + \lb_1 \, \Zc (s) + 
\rho (s)
\ee
where the $\lb_j$ are scalars and $\rho$ is a holomorphic function of $s 
\in \Cb$, thus proving theorem 2.
\medskip

\section{The Local Index Formula for $SU_q(2)$, ($q=0$).}
\smallskip

\noindent In this section we shall compute the local index formula 
for the above spectral triple. Since the dimension spectrum is 
simple and equal to $\{ 1,2,3 \} \sbs 
\Cb$ the cyclic cocycle given by the local index formula has 
two components $\vp_1$ and $\vp_3$ of degree $1$ and $3$ given,
 up to an overall multiplication by $( 2 i \pi)^{1/2}$, by
\begin{eqnarray}
\label{eq1.70}
\vp_1 (a^0 , a^1) = \int\!\!\!\!\!\!- \ a^0 [D ,a^1 ]\vert D 
\vert^{-1} - \frac{1}{4} \int\!\!\!\!\!\!- \ a^0 \nabla ([D ,a^1 ]) 
\, \vert D \vert^{-3} \\
+ \frac{1}{8} \int\!\!\!\!\!\!- \ a^0 \nabla^2 
([D ,a^1 ]) \, \vert D \vert^{-5} \nonumber 
\end{eqnarray}
and,
\be
\label{eq1.71}
\vp_3 (a^0 , a^1 , a^2 , a^3) = \frac{1}{12} 
\int\!\!\!\!\!\!- \ a^0 [D ,a^1 ][D , a^2] [D , a^3 ]\, \vert 
D \vert^{-3} \, .
\ee
With these notations the cocycle equation is,
\be
\label{eq1.72}
b \vp_1 + B \vp_3 = 0 \, .
\ee
The following formulas define a cyclic cocycle $\tau_1$ on $\Ac$,
\be
\label{eq1.49}
\tau_1 (\a^{*k}, x)=\tau_1 (x,\a^{*k})= \tau_1 (\a^{l}, x)=\tau_1 (x,\a^{l})
=0, 
\ee
for all integers $k$, $l$ and any $x \in \Ac$,
\be
\tau_1 (\a^{*k} f(\b) \a^{\ell} , \a^{*k'}  g(\b) \a^{\ell'}) = 0
\ee
unless  $\ell'=k$, $k'=\ell$ and
$$
\tau_1 (\a^{*k} f(\b)\a^{\ell} , \a^{*\ell} g(\b) \a^k) = \frac{1}{\pi i} \int_{S^1} f \,\,
 {\rm d} g \, 
 \, .
$$
Let $\vp_0$ be the 0-cochain given by $\vp_0 (\a^{*k} f(\b) \a^{\ell}) = 0 
 \,\, \hbox{unless} \,\, k = \ell $ and,
\be
\label{eq1.89}
\vp_0 (\a^{*k} f(\b)
\a^k) =  \rho (k)\, \, \frac{1}{2\pi}  \int_{S^1} f \, \,{\rm d} \t ,
\ee
where $\rho (j) =\frac{2}{3} - j - j^2$.
Finally, let $\vp_2$ be the 2-cochain given by the pull back by $\s$ of the 
cochain  $\frac{-1}{24}\frac{1}{2\pi i} \, \int f_0 
f'_1 f''_2 \, {\rm d} \t $ on  $C^{\ify} (S^1)$.

\noindent Our next task is to prove the following result,
\bigskip

\noindent {\bf Theorem 3.} {\it The local index formula of the spectral 
triple $(\Ac , \Hc , D)$ is given by the cyclic cocycle $\tau_1$ up to the coboundary of the cochain ($\vp_0$,$\vp_2$}).

\bigskip

\noindent The precise equations are,
\be
\label{eq1.272}
\vp_1= \,\tau_1 + b \vp_0 + B \vp_2 \,, \qquad \vp_3= \, b \vp_2 \,.
\ee
\noindent The proof is a computation but we shall go through it in details in order to get familiar with various ways of 
computing residues and manipulating "infinitesimals"
in the sense of the quantized calculus. In other words our purpose is not concision but rather a leisurly account of the details.

\subsection{Restriction to $C^{\ify}(\b)$}

\noindent Let us first concentrate on the restriction of the cocycle $\vp$
to the subalgebra $C^{\ify}(\b)$ generated by $\b$ and  $\b^{*}$. To see the subspace of $\Hc$ responsible 
for the non-triviality of that cocycle we follow the action of $\b$ on the vectors,
\be 
\label{eq1.7}
\xi_{-n} = e_{-\frac{n}{2} , \frac{n}{2}}^{\left( \frac{n}{2} 
\right)} \qquad n \geq 0 \, , \ n \in \Nb \, .
\ee
and,
\be
\label{eq1.8}
\xi_n = e_{n/2 , -n/2}^{(n/2)} \qquad n \geq 0 \, , \ n \in \Nb
\ee

\noindent For $n>0$, (\ref{eq21.3}) shows that $\b (\xi_{-n}) = \xi_{-(n-1)} = 
\xi_{-n+1}$, with $\b (\xi_{-1}) = e_{0,0}^{(0)} = \xi_0$.
Next, $\b (\xi_0) = \b (e_{0,0}^{(0)}) = -e_{(1/2) , - 
1/2}^{(1/2)}= -\xi_1$, and for $n>0$ (\ref{eq21.4}) shows that
$\b (\xi_n) = -\xi_{n+1}$. Thus,
\be
\label{eq1.9}
\b (\xi_n) = - {\rm sign} (n)\, \, \xi_{n+1} \qquad ({\rm sign} (0) = 
1)
\ee
We let  $\ell^2 (\Zb) = 
\Hc_0 \sbs \Hc$ be the subspace of $\Hc$ spanned by the $ \xi_n$
and rewrite the above equality as,
\be
\label{eq1.10}
\b = -UH \ {\rm on} \ \ell^2 (\Zb) \sbs \Hc
\ee
where $H$ is the sign operator and $U$ the shift,
\be
\label{eq1.11}
U \xi_n = \xi_{n+1} \, .
\ee
The operator $D$ also restricts to the subspace $\ell^2 (\Zb) = 
\Hc_0 \sbs \Hc$ and its restriction $D_0$ is given by,
\be
\label{eq1.12}
D_0 \xi_n = {\rm sign} (n) \, \vert n \vert \, \xi_n = 
n \, \xi_n \qquad \fl n \, .
\ee
The unitary $W=  \,e^{i \frac{\pi}{2}( \vert D_0 \vert - D_0)}$
commutes with $D_0$ and conjugates $U$ to $-UH$,
\be
\label{eq1.16}
W \,U \, W^{*}=\, -UH 
\ee
Thus the triple ($\b$, $\Hc_0$,$D_0$)  is isomorphic to,
\be
\label{eq1.13}
\left( e^{i\t} , \, \,L^2 (S^1) ,  \, (-i) \, 
\frac{\partial}{\partial \,\t} \right) \, .
\ee

\noindent In particular the index and cyclic cohomology pairings with the 
restriction to $\Hc_0$ are non trivial and we control,
\be
\label{eq1.18}
{\rm Res}_{s=1} \, {\rm Trace}_{\Hc_0} \, (\b^* [D_0,\b] \,  \vert D_0\vert^{-s}) \,=2 .
\ee
This however does not suffice to get the non-triviality of the restriction of $\vp$
to $C^{\ify}(\b)$ since we need to control the residues on $e \Hc$ where,
as above $e$ is the support of $\b$.  To see what happens
we shall conjugate the restriction of both $\b$ and $D$ to
the orthogonal complement of $\Hc_0$ in $e \Hc$ with a very
simple triple.
Let us define for each $k \in \Nb$ the vectors,
\be
\label{eq1.19}
\xi_{-n}^{(k)} = e_{-\left( \frac{k+n}{2} \right) , 
\frac{n-k}{2}}^{\left( \frac{k+n}{2}\right)} \qquad n \geq 0
\ee
$$
\xi_n^{(k)} = e_{\frac{n-k}{2} , - \left( \frac{k+n}{2} 
\right)}^{\frac{k+n}{2}} \qquad n \geq 0
$$
so that $\xi_0^{(k)} = e_{-k/2 , -k/2}^{(k/2)}$.

\smallskip

\noindent For $n > 0$ one has
$$
\b (\xi_{-n}^{(k)}) = \b \left( e_{- \left( \frac{k+n}{2} 
\right) , \frac{n-k}{2}}^{\left( \frac{k+n}{2}\right)} \right) = 
e_{- \left( \frac{k+n-1}{2} \right) , \frac{n-1-k}{2}}^{\left( 
\frac{k+n-1}{2}\right)} = \xi_{-(n-1)}^{(k)} = \xi_{-n+1}^{(k)} 
\, .
$$
For $n=0$,
$$
\b (\xi_0^{(k)}) = \b (e_{-k/2 , -k/2}^{(k/2)}) = 
-e_{-\frac{k}{2} + \frac{1}{2} , - \left( \frac{k+1}{2} 
\right)}^{\left( \frac{k+1}{2} \right)} = \xi_1^{(k)}
$$
and,
$$
\b (\xi_n^{(k)}) = \b \left( e_{\frac{n-k}{2} , - \left( 
\frac{n+k}{2} \right)}^{\left( \frac{k+n}{2}\right)} \right) =  
- e_{\frac{n+1-k}{2} , - \left( \frac{n+1+k}{2} \right)}^{\left( 
\frac{k+n+1}{2}\right)} = -\xi_{n+1}^{(k)} \, .
$$
Thus, as in (\ref{eq1.9}) we have,
\be
\label{eq1.20}
\b (\xi_n^{(k)}) = - {\rm sign} \, (n) \, \xi_{n+1}^{(k)} \, .
\ee
Now $\b (e_{ij}^{(m)}) = 0$ unless $i = -m$ or $j = -m$ and for 
any $m \in \frac{1}{2} \, \Nb$ the vectors $e_{-m,j}^{(m)}$ and 
$e_{i,(-m)}^{(m)}$ are of the form $\xi_n^{(k)}$. Indeed in the 
first case one takes $n=m+j$, $k = m-j$ which are both in $\Nb$, 
and $\xi_{-n}^{(k)} = e_{-m,j}^{(m)}$. In the second case 
$n=m+i$, $k=m-i$ are both in $\Nb$ and $\xi_n^{(k)} = 
e_{i,-m}^{(m)}$.

\smallskip

\noindent We then let $\Hc_k$ be the span of the $\xi_n^{(k)}$, $n 
\in \Zb$ and
\be
\label{eq1.21}
\Hc' = \bigoplus_{k \geq 1} \Hc_k = \ell^2 (\Zb) \otimes \ell^2 
(\Nb^+) \, .
\ee
The operator $D$ restricts to $\Hc'$ and is given there by,
\be
\label{eq1.22}
D'=\vert D_0 \vert \otimes 1 + 1 \otimes N
\ee
where $N$ is the number operator $N \ve_k = k \ve_k$.

\smallskip

\noindent Also $\b$ is $-UH \otimes 1$ and we can conjugate it 
as in (\ref{eq1.16})  back to $U \otimes 1$.

\smallskip

\noindent Thus the triple $(\b,\Hc', D')$ is isomorphic to 
\be
\label{factor}
(U \otimes 1,\,\,\ell^2 (\Zb) \otimes \ell^2 (\Nb^{+}),\,\,\vert D_0 \vert \otimes 1 + 1 \otimes N) 
\ee
\smallskip

\noindent The metric dimension is 2 in this situation,
 and the contribution $\vp'$  of $\Hc'$ to the restriction of $\vp$ to $C^{\ify}(\b)$ only has 
a one dimensional component $\vp'_1$ which involves the 
two terms,
\be
\label{eqvp}
\vp'_1 (a^0 , a^1) = \int\!\!\!\!\!\!- \ a^0 [D' ,a^1 ] \vert D' 
\vert^{-1} - \frac{1}{4} \int\!\!\!\!\!\!- \ a^0 \nabla ([D' ,a^1 ]) 
\, \vert D' \vert^{-3} ,\,\,\,a^j \in C^{\ify}(\b)
\ee
Since $D'$ is positive, it is K-homologically trivial and the
above cocycle must vanish identically on $C^{\ify}(\b)$. 
As we shall see this vanishing
holds because of the precise ratio $-\frac{1}{4}$
of the coefficients in the local index formula.

\noindent To see this, we need to compute the poles and residues
of functions of the form ${\rm Trace} \, ((T \otimes 1) \, \vert D' \vert^{-s})$
for operators $T$ in $\ell^2 (\Zb)$. 
 For that purpose it is most efficient to use the
well known relation between residues of zeta functions and asymptotic
expansions  of related theta functions. More specifically,
 for $\lb >0$ and $Re(s)>0$, the equality,
\be
\label{eq1.23}
\lb^{-s} = \frac{1}{\G (s)} \int_0^{\ify} e^{-t\lb} \, t^s \, 
\frac{dt}{t}
\ee
gives
\begin{eqnarray}
{\rm Trace} \, ((T \otimes 1) \, \vert D' \vert^{-s}) &= 
&\frac{1}{\G (s)} \int_0^{\ify} {\rm Trace} \, ((T \otimes 1) \, 
e^{-tD'}) \, t^s \, \frac{dt}{t} \nonumber \\
&= &\frac{1}{\G (s)} \int_0^{\ify} {\rm Trace} \, (T e^{-t \vert 
D_0 \vert}) \left( \frac{1}{e^t -1} \right) t^s \, \frac{dt}{t} \, 
. \nonumber
\end{eqnarray}
Thus if we assume that one has an expansion of the form ${\rm 
Trace} \, (T e^{-t \vert D_0 \vert}) = \frac{a}{t} + b + c  \,t + 0(t^2)$ one 
gets, using
\be
\label{eq1.24}
\frac{1}{e^t - 1} = \frac{1}{t} - \frac{1}{2} + \frac{t}{12}+ 0(t^2)
\ee
the equality modulo holomorphic functions of $s$, $Re(s)>0$,
$$
{\rm Trace} \, ((T \otimes 1) \, \vert D' \vert^{-s})
\sim \frac{1}{\G (s)} \int_0^1 \vp (t) \, t^s \, \frac{dt}{t} ,
$$
where,
$$
\vp(t) = \frac{a}{t^2} + \frac{\left( b - 
\frac{a}{2} \right)}{t} + \frac{a}{12}-\frac{b}{2}+c  \, .
$$
One has $\int_0^1  t^{\a} \, \frac{dt}{t} = \frac{1}{\a}$, thus 
one gets 2 poles $s=2$ and $s=1$, and the expansion,
$$
{\rm Trace} \, ((T \otimes 1) \, \vert D' \vert^{-s}) \sim 
\frac{a}{\G (s)(s-2)} + \left( b - \frac{a}{2} \right) 
\frac{1}{\G (s) (s-1)} + \cdots
$$
so that, 
\begin{eqnarray}
\label{eq1.25}
{\rm Res}_{s=2} ({\rm Trace} \, ((T \otimes 1) \vert D' 
\vert^{-s})) = a \\
{\rm Res}_{s=1} ({\rm Trace} \, ((T \otimes 1) \vert D' 
\vert^{-s})) =  b - \frac{a}{2} 
\end{eqnarray}
Let us compute $\vp'_1(\b^*, \b)$. The first term in (\ref{eqvp})
is $\int\!\!\!\!\!\!- \ \b^* [D' ,\b ] \vert D' 
\vert^{-1}$. Thus we take $T = U^* [\vert D_0 \vert,U]$. 
One has $[\vert D_0 \vert,U]=U H$, $T=H$ and
\be
\label{eta}
{\rm Trace} \, (T e^{-t \vert D_0 \vert}) = \sum_{\Zb} {\rm sign} 
\, (n) \, e^{-t \vert n \vert} = - \sum_1^{\ify} e^{-tk} + 1 + 
\sum_1^{\ify} e^{-tk} = 1 \, .
\ee
Thus in that case $a=0$, $b=1$ and,
\be
\label{eq1.26}
\int\!\!\!\!\!\!- \ \b^* [D' ,\b ] \vert D' 
\vert^{-1} = 1 \, .
\ee
The second term in (\ref{eqvp}) is 
$\int\!\!\!\!\!\!- \ \b^*  \nabla ([D' ,\b  ]) 
\, \vert D' \vert^{-3}$ where $\nabla $ is the 
commutator with $ D'^2$. If we let as above $\d$
be the commutator with $\vert D' \vert$, one has
$$
\nabla T =  
\d (T) \, \vert D' \vert + \vert D' \vert \, \d (T) 
\, .
$$
Thus, permuting  $\vert D' \vert$ modulo operators
of lower order, we get,
\be
\label{eqdouble}
\int\!\!\!\!\!\!- \ \b^*  \nabla ([D' ,\b  ]) 
\, \vert D' \vert^{-3}=\,2 \,\
\int\!\!\!\!\!\!- \ \b^*  \d ([D' ,\b  ]) 
\, \vert D' \vert^{-2}
\ee
To compute the r. h. s. we take $T = U^* \d^2 (U)$ and look at the residue at 
$s=2$. One has $\d (U) = UH$, $\d^2 (U) = (UH)H = U$ and $T=1$. 
Thus 
\begin{eqnarray}
{\rm Trace} \, (T e^{-t \vert D_0 \vert}) &= &\sum_{\Zb} e^{-t 
\vert n \vert} = 1 + 2 \sum_1^{\ify} e^{-tn} \nonumber \\
&= & 1+ \frac{2}{e^t 
-1} \sim \frac{2}{t} + 0(t) \nonumber
\end{eqnarray}
so that $a=2$, $b=0$. Thus we get,
\be
\label{eq1.28}
\int\!\!\!\!\!\!- \ \b^*  \d ([D' ,\b  ]) 
\, \vert D' \vert^{-2}
={\rm Res}_{s=2} ({\rm Trace} \, ((U^* \d^2 (U)\otimes 1) \, \vert D' 
\vert^{-s}) = 2
\ee
Thus we get,
\be
\label{eq1.34'}
\int_{\Hc'}\!\!\!\!\!\!\!\!\!\!- \ \ \b^*\,[D,  \b ]\,\vert D \vert^{-1}
 = 1 \, ,
\
\int_{\Hc'}\!\!\!\!\!\!\!\!\!\!- \ \ \b^*\,\nabla ([D,  \b ])\,\vert D \vert^{-3}
 = 4 
\ee
and $\vp'_1(\b^*,\b)=0$ precisely because of the coefficient $\frac{-1}{4}$
in (\ref{eqvp}).

\medskip

\noindent One proceeds similarly to compute $\vp'_1(\b, \b^*
)$. The first term in (\ref{eqvp}) comes from $T = U \d (U^*)$.
 In the canonical 
basis $\ve_n $ of $\ell^2 (\Zb)$ one has $U \d 
(U^*) \, \ve_n = (\vert n-1 \vert - \vert n \vert) \, \ve_n
= - {\rm sign} \, (n-1) \, \ve_n$
and ${\rm Trace} \, (T e^{-t \vert D_0 \vert}) = \sum (-{\rm sign} 
\, (n-1)) \, e^{-t \vert n \vert} = 1$. 
Thus,
\be
\int\!\!\!\!\!\!- \ \b \,\,[D' ,\b^* ] \vert D' 
\vert^{-1} = 1 \, .
\ee

\smallskip

\noindent Since $U \d^2 (U^*)=1$ the computation of the second term of 
(\ref{eqvp}) is the same as above and we get,
\be
\label{eq1.34''}
\int_{\Hc'}\!\!\!\!\!\!\!\!\!\!- \ \ \b\,[D,  \b^* ]\,\vert D \vert^{-1}
 = 1 \, ,
\
\int_{\Hc'}\!\!\!\!\!\!\!\!\!\!- \ \ \b\,\nabla ([D,  \b^* ])\,\vert D \vert^{-3}
 = 4 
\ee
so that $\vp'_1(\b,\b^*)=0$.

\noindent Let us now take $n > 
0$ and compute $\vp'_1(\b^{*n},\b^n)$.
The first term of (\ref{eqvp}) involves $T=U^{*n} \d (U^n)$. One has, 
$$
\d (U^n) = UHU^{n-1} + U^2 HU^{n-2} + \cdots + U^j H U^{n-j} + 
\cdots + U^n H = \sum_{j=1}^n U^j H U^{n-j} \, ,
$$
$$
U^{*n} \d (U^n) = (U^*)^{n-1} H U^{n-1} + \cdots + H \, .
$$
One has $U^{*k} H U^k \ve_{\ell} = {\rm sign} (k+\ell) \, 
\ve_{\ell}$ and,
$$
{\rm Trace} \, (U^{*k} H U^k e^{-t\vert D_0 \vert}) = {\rm Trace} 
\, (H e^{-t\vert D_0 \vert}) + 2 \sum_1^k e^{-t \vert j \vert} \sim 
1+2k+0(t) \, .
$$
Thus,
$$
{\rm Trace} \, (U^{*n} \d (U^n) \, e^{-t \vert D_0 \vert}) \sim 
\sum_{j=1}^n (2(n-j)+1)+0(t) = \sum_{k=0}^{n-1} (2k+1) + 0(t) = 
n^2 + 0(t) \, ,
$$
and $(n > 0)$
\be
\label{eq1.36}
\int\!\!\!\!\!\!- \ \ \b^{*n}\,[D',  \b^n ]\,\vert D' \vert^{-1}
 = n^2 \, ,
\ee
Now modulo finite rank operators one has $\d (U^n) = n U^n H$ and 
$\d^2 (U^n) = n^2 U^n$, thus as above,
\be
\label{eq1.37}
\int\!\!\!\!\!\!- \ \ \b^{*n}\,\nabla ([D',  \b^n ])\,\vert D' \vert^{-3}
 = 4 n^2,
\ee
so that $\vp'_1(\b^{*n},\b^n)=0$.

\medskip
\noindent Finally the computation of $\vp'_1(\b^n,\b^{*n})$ involves $T=U^n \d (U^{*n})$. One has
$$
U^n \d (U^{*n}) = - \sum_{k=1}^n U^k H U^{*k} \, ,
$$ and $U^k H U^{*k} \ve_{\ell} = {\rm sign} \, (\ell - k) \, 
\ve_{\ell}$
so that,
$$
{\rm Trace} \, (U^k H U^{*k} e^{-t \vert D_0 \vert}) - {\rm Trace} \, 
(H e^{-t \vert D_0 \vert}) = -2 \sum_0^{k-1} e^{-tj}
$$
and,
$$
{\rm Trace} \, (U^k H U^{*k} e^{-t \vert D_0 \vert}) \sim 1-2k + 0(t) 
\, .
$$
Thus, 
\begin{eqnarray}
{\rm Trace} \, (U^n \d (U^{*n}) \, e^{-t \vert D_0 \vert}) &= &- 
\sum_{k=1}^n {\rm Trace} \, (U^k H U^{*k} e^{-t \vert D_0 \vert}) 
\nonumber \\
&\sim &- \sum_{k=1}^n (1-2k) + 0(t) \sim n^2 + 0(t) \, , \nonumber
\end{eqnarray}
and,
\be
\label{eq1.38}
\int\!\!\!\!\!\!- \ \ \b^n\,[D',  \b^{*n} ]\,\vert D' \vert^{-1}
 = n^2 \, .
\ee
Also as above,
\be
\label{eq1.39'}
\int\!\!\!\!\!\!-  \ \ \b^n\,\nabla ([D',  \b^{*n} ])\,\vert D' \vert^{-3}
 = 4 n^2.
\ee
so that we get the required vanishing,
$$
\vp'_1=0
$$
What is instructive in the above computation is that this vanishing
which is required by theorem 1, involves because of the factorisation
(\ref{factor}) terms such as "$ {\rm Trace}(H)$" which appear in equation (\ref{eta}) and are similar to eta-invariants.

\noindent We have thus shown that $\vp_1=\tau_1$ on $C^{\ify}(\b)$, or equivalently that,
\be
\label{eqphi}
\int\!\!\!\!\!\!-  \ \ \b^{-n}\,[D,  \b^n ]\,\vert D \vert^{-1}- \frac{1}{4}\int\!\!\!\!\!\!-  \ \ \b^{-n}\,\nabla ([D,  \b^n ])\,\vert D \vert^{-3}
 = 2n.
\ee

\medskip
\subsection{Restriction to the ideal $\Jc$}

Let us extend this computation to $\Jc = 
{\rm Ker} \, \s$. The component $\vp_3$ vanishes on  $\Jc$ and we
just need to compute $\vp_1=\vp_1^{(0)} - \frac{1}{4} \, \vp_1^{(1)} + \frac{1}{8} \, 
\vp_1^{(2)}$. We begin by $\vp_1^{(0)}(\mu' , \mu)=\int\!\!\!\!\!\!- \ \mu' [D, \mu] \vert D 
\vert^{-1} $, and need only consider the case where $\mu = \a^{*k} \b^{n} \a^{\ell}$
and $\mu' = \a^{*k'} \b^{n'} \a^{\ell'}$ are monomials. (As above $\b^{-n}=\b^{*n}$ and
  $\b^0=e$).

\noindent  With $F=\,2P-1$ one has $[F,\a] = 0$ and
\be
\label{eq1.74}
[D,\a^{\ell}] = - \ell F \a^{\ell} = - \ell \a^{\ell} F
\ee
and
\be
\label{eq1.75}
[D,\a^{*k}] = k F \a^{*k} = k \a^{*k} F \, .
\ee
Thus,
\begin{eqnarray}
[D,\mu] &= &[D,\a^{*k}] \, \b^n \a^{\ell} + \a^{*k} [D,\b^n] \, 
\a^{\ell} + \a^{*k} \b^n [D,\a^{\ell}] \nonumber \\
&= &k \a^{*k} F \b^n \a^{\ell} + \a^{*k} [D,\b^n] \, \a^{\ell} - \ell 
\a^{*k} \b^n F \a^{\ell} \, , \nonumber
\end{eqnarray}
so that
\be
\label{eq1.76}
[D,\mu] = \a^{*k} (k F \b^n + [D,\b^n] - \ell \b^n F) \, \a^{\ell}
\ee

\smallskip

\noindent The bigrading (\ref{eq1.43}) shows that 
$\int\!\!\!\!\!- \ \mu' [D,\mu] \, \vert D \vert^{-1}$ vanishes unless 
both total degrees are 0, i.e.
\be
\label{eqhom}
\ell + \ell' - k - k' =0 \, , \ n' = -n \, .
\ee
The element $X=(k F \b^n + [D,\b^n] - \ell \b^n F)$ satisfies
$eXe=X$ so that the product 
$$
\mu' [D,\mu] = \a^{*k'} \b^{n'} \a^{\ell'} \a^{*k}\,X\, \a^{\ell},
$$
vanishes unless $\ell' = k$. Combining with (\ref{eqhom}) we 
get $\ell = k'$, and can assume that $\mu'=\a^{*\ell} \b^{*n} \a^k $.
Then $\mu'[D,\mu]=\a^{*\ell} \b^{*n}\,X  \, \a^{\ell}$
so that we just need to compute,  $$
\int\!\!\!\!\!\!- \ \a^{*\ell} \b^{*n}\, X \, \a^{\ell} \, \vert D 
\vert^{-1} \, .
$$
Now by (\ref{eqperm}) one has,

\be
\label{eq1.79}
\vert D \vert^{-1} \a^{*\ell} = \a^{*\ell} \vert D \vert^{-1} - \ell 
\a^{*\ell} \vert D \vert^{-2} + 0 (\vert D \vert^{-3})
\ee
Thus,
\begin{eqnarray}
\label{eqmo}
&&\int\!\!\!\!\!- \ \mu' [D,\mu] \, \vert D \vert^{-1}
=\int\!\!\!\!\!\!- \ \b^{*n} X \, \vert D \vert^{-1} - \ell 
\int\!\!\!\!\!\!- \ \b^{*n} X \, \vert D \vert^{-2}
\end{eqnarray}
where
$$
X = (k F \b^n + [D,\b^n] - \ell \b^n F) \, .
$$
Note that $\a , \a^*$ have now disappeared so that we can compute 
using the subspaces $\Hc_0$ and $\Hc'$ of $e\Hc$. Note also that on 
$\Hc'$ one has $F=-1$ since $P=0$. Only $\Hc'$ matters for 
$\int\!\!\!\!\!- \ \b^{*n} X \, \vert D \vert^{-2}$. One has 
$\int\!\!\!\!\!- \ \b^{*n} k F \b^n \vert D \vert^{-2} = -k 
\int_{\Hc'}\!\!\!\!\!\!\!\!\!\!- \ \ \vert D \vert^{-2}$ since $F \sim 
-1$, and $\int\!\!\!\!\!- \ \b^{*n} (-\ell \b^n F) \, \vert D 
\vert^{-2} = \ell \int_{\Hc'}\!\!\!\!\!\!\!\!\!\!- \ \ \vert D 
\vert^{-2}$. Also $\int\!\!\!\!\!- \ \b^{*n} [D,\b^n] \, \vert D 
\vert^{-2} = 0$ since $\int_{L^2 
(S^1)}\!\!\!\!\!\!\!\!\!\!\!\!\!\!\!\!\!\!\!- \qquad \ U^{-n} [\vert D_0 
\vert , U^n] \, \vert D_0 \vert^{-1} = 0$. Thus,
\be
\label{eq1.80}
 \int\!\!\!\!\!\!- \ \b^{*n} X \, \vert D \vert^{-2} =\,
(\ell - k) \int_{\Hc'}\!\!\!\!\!\!\!\!\!\!- \ \ \vert D \vert^{-2} \, 
.
\ee
One has
$$
\int\!\!\!\!\!\!- \ \b^{*n} X \, \vert D \vert^{-1} = 
\int\!\!\!\!\!\!- \ \b^{*n} k F \b^n \vert D \vert^{-1} + 
\int\!\!\!\!\!\!- \ \b^{*n} [D,\b^n] \, \vert D \vert^{-1} - \ell 
\int\!\!\!\!\!\!- \ e F \, \vert D \vert^{-1}
$$
where the $\int\!\!\!\!\!-$ are on $\Hc' \oplus \Hc_0$. One has 
$\int_{\Hc_0}\!\!\!\!\!\!\!\!\!\!\!- \ \ \ \ \b^{*n} F \b^n \vert D 
\vert^{-1} = 0$ and $\int_{\Hc_0}\!\!\!\!\!\!\!\!\!\!\!- \ \ \ F \, \vert 
D \vert^{-1} = 0$. On $\Hc'$ one has $F=-1$, thus,
\be
\label{eq1.81}
\int\!\!\!\!\!\!- \ k \b^{*n} F \b^n \vert D \vert^{-1} = -k 
\int_{\Hc'}\!\!\!\!\!\!\!\!\!\!- \ \ \vert D \vert^{-1}
\ee
and
\be
\label{eq1.82}
-   \ell  \int\!\!\!\!\!\!- \ e F \, \vert D \vert^{-1} = \ell 
\int_{\Hc'}\!\!\!\!\!\!\!\!\!\!- \ \ \vert D \vert^{-1}
\ee
so that,
\be
\label{eq1.83}
\int\!\!\!\!\!\!- \ \b^{*n} X \, \vert D \vert^{-1} = 
\int\!\!\!\!\!\!- \ \b^{*n} [D,\b^n] \, \vert D \vert^{-1} + (\ell - 
k) \int_{\Hc'}\!\!\!\!\!\!\!\!\!\!- \ \ \vert D \vert^{-1} \, .
\ee
We now need to compute $- \frac{1}{4} \ \int\!\!\!\!\!\!- \ \mu' 
\nabla ([D,\mu]) \, \vert D \vert^{-3}$ with $\mu , \mu'$ monomials in 
$\Jc$ as above. 

\noindent Since $\vert D \vert^{-2}$ is order 1 on $\Jc$ we can 
replace the above by $- \frac{1}{2} \ \int\!\!\!\!\!\!- \ \mu' \d 
([D,\mu])$ $\vert D \vert^{-2}$ where $\d (x) = [\vert D \vert , x ]$. 
Moreover, using (\ref{eq1.76}),
$$
\d ([D,\mu]) = \d (\a^{*k} X \a^{\ell}) = k \a^{*k} X \a^{\ell} + 
\a^{*k} \d (X) \, \a^{\ell} - \ell \a^{*k} X \a^{\ell} \, ,
$$
\be
\label{eq1.84}
\d ([D,\mu]) = (k-\ell) \, \a^{*k} X \a^{\ell} + \a^{*k} \d (X) \, 
\a^{\ell}
\ee
with $X = k F \b^n + [D,\b^n] - \ell \b^n F$. 

\noindent As above for 
$\vp_1^{(0)}$ we get that $\vp_1^{(1)}$ vanishes unless $\ell' = k$, 
$k' = \ell$, $n' = -n$ so that $\mu' = \a^{*\ell} \b^{*n} \a^k$ and we 
can replace $\mu' \d ([D,\mu]) \, \vert D \vert^{-2}$ by
\be
\label{eq1.85}
(k-\ell) \, \b^{*n} X \, \vert D \vert^{-2} + \b^{*n} \d (X) \, \vert 
D \vert^{-2} \, .
\ee
Now by (\ref{eq1.80}), 
$$
\int\!\!\!\!\!\!- \ \b^{*n} X \, \vert D \vert^{-2} = (\ell-k) 
\int_{\Hc'}\!\!\!\!\!\!\!\!\!\!- \ \ \vert D \vert^{-2} \, .
$$
Moreover one has $\int\!\!\!\!\!- \ \b^{*n} \d (F \b^n) \, \vert D 
\vert^{-2}=\int\!\!\!\!\!- \ \b^{*n} \d (\b^n F) \, \vert D 
\vert^{-2} = 0$ since only the $\int\!\!\!\!\!-$ on $\Hc'$ matters 
and $\int_{L^2 
(S^1)}\!\!\!\!\!\!\!\!\!\!\!\!\!\!\!\!\!\!\!- \qquad \ U^{-n} [\vert D_0 
\vert , U^n] \, \vert D_0 \vert^{-1} = 0$. Thus,
\begin{eqnarray}
\label{eq1.86}
- \frac{1}{4} \ \int\!\!\!\!\!\!- \ \mu' \nabla ([D,\mu]) \, \vert D 
\vert^{-3} &= &- \frac{1}{2} \, (k-\ell)(\ell-k) 
\int_{\Hc'}\!\!\!\!\!\!\!\!\!\!- \ \ \vert D \vert^{-2} \nonumber \\ &- &\frac{1}{2} \ \int\!\!\!\!\!\!- \ \b^{*n} \d ([D,\b^n]) \, \vert D 
\vert^{-2} \, .
\end{eqnarray}
Since $\vp_1^{(2)}=0$ on $\Jc$ we get,
\begin{eqnarray}
\vp_1 (\mu',\mu) &= &\int\!\!\!\!\!\!- \ \b^{*n} [D,\b^n] \, \vert D 
\vert^{-1} + (\ell-k) \int_{\Hc'}\!\!\!\!\!\!\!\!\!\!- \ \ \vert D 
\vert^{-1} - \ell \,\,(\ell-k) \int_{\Hc'}\!\!\!\!\!\!\!\!\!\!- \ \ \vert 
D \vert^{-2} \nonumber \\
& - &\frac{1}{2} \, (k-\ell)(\ell-k) \int_{\Hc'}\!\!\!\!\!\!\!\!\!\!- 
\ \ \vert D \vert^{-2} - \frac{1}{2} \ \int\!\!\!\!\!\!- \ \b^{*n} \d 
([D,\b^n]) \, \vert D \vert^{-2} \, . \nonumber
\end{eqnarray}
Now by (\ref{eqphi}),
\be
\label{eq1.87}
\int\!\!\!\!\!\!- \ \b^{*n} [D,\b^n] \, \vert D \vert^{-1} 
-\frac{1}{2} \ \int\!\!\!\!\!\!- \ \b^{*n} \d ([D,\b^n]) \, \vert D 
\vert^{-2} = 2n
\ee
 Thus we get,
\be
\label{eq1.88}
\vp_1 (\mu',\mu) = 2n + (\ell-k) \int_{\Hc'}\!\!\!\!\!\!\!\!\!\!- \ \ 
\vert D \vert^{-1} - \frac{1}{2} \, (\ell^2 - k^2) 
\int_{\Hc'}\!\!\!\!\!\!\!\!\!\!- \ \ \vert D \vert^{-2} \, .
\ee
Let us show that $\vp_1$ is cohomologous to $\tau_1$ on $\Jc$. 
Indeed, let $\rho (k)$ be an arbitrary sequence of
polynomial growth and $\vp_0$ be the 
0-cochain given by,
\be
\label{eq1.89}
\vp_0 ( \a^{*k} \b^n \a^{\ell}) = 0 
\, \, \, \hbox{unless} \, \, k = \ell , \qquad \vp_0 (\a^{*k} \b^0 
\a^k) = \rho (k) \, .
\ee
Then 
$$
(b \vp_0)(\mu',\mu) = \vp_0 (\mu'\mu) - \vp_0 (\mu\mu')
$$
and both terms vanish unless $k=\ell'$, $k'=\ell$, $n'=-n$. Moreover in 
that case
$$
\mu\mu' = \a^{*k} \b^n \a^{\ell} \a^{*k'} \b^{n'} \a^{\ell'} = \a^{*k} 
\b^0 \a^k
$$
while
$$
\mu'\mu = \a^{*\ell} \b^{-n} \a^k \a^{*k} \b^n \a^{\ell} = \a^{*\ell} 
\b^0 \a^{\ell}
$$
so that
$$
(b \vp_0)(\mu' , \mu) = \rho (\ell) - \rho (k) \, . $$
Thus, with,
\be
\label{eq1.90}
\rho (k) = k \int_{\Hc'}\!\!\!\!\!\!\!\!\!\!- \ \ \vert D \vert^{-1} - 
\frac{1}{2} \, k^2 \int_{\Hc'}\!\!\!\!\!\!\!\!\!\!- \ \ \vert D 
\vert^{-2}
\ee
we have, on the ideal $\Jc$,
\be
\label{eq1.91}
\vp_1 = \tau_1 + b \vp_0 \, .
\ee
Let us now extend this equality to the case when only one of the 
variables $\mu,\mu'$ belongs to the ideal $\Jc$.

\smallskip

\noindent Assuming first that $\mu$ belongs to the ideal $\Jc$,
we just need to compute  $\vp_1 (\mu',\mu)$  for $\mu = \a^{*k} \b^0 
\a^{\ell}$ and $\mu' = \a^{*k'}$ if $k' = \ell - k \geq 0$ or $\mu' = 
\a^{\ell'}$ if $\ell' = k - \ell > 0$. One has by  (\ref{eq1.76}),
$[D,\mu] = \a^{*k} X 
\a^{\ell}$, $X = (k-\ell) \, \b^0 F$ since $[D,\b^0] = 0$ and 
$[F,\b^0] = 0$. Thus, for $k' \geq 0$,
$$
\int\!\!\!\!\!\!- \ \mu' [D,\mu] \, \vert D \vert^{-1} =  \int\!\!\!\!\!\!- \ \a^{*k'} \a^{*k} (k-\ell) \, \b^0 F \a^{\ell} \, 
\vert D \vert^{-1} = (k-\ell) \int\!\!\!\!\!\!- \ \a^{*\ell} \b^0 F 
\a^{\ell} \, \vert D \vert^{-1}
$$
since $k+k'=\ell$.

\smallskip

\noindent For $k' < 0$, $\ell' = k-\ell > 0$ one gets 
$$
\int\!\!\!\!\!\!- \ \a^{\ell'} \a^{*k} (k-\ell) \b^0 F \a^{\ell} \, 
\vert D \vert^{-1} = (k-\ell) \int\!\!\!\!\!\!- \ \a^{*\ell} \b^0 F 
\a^{\ell} \, \vert D \vert^{-1} \, .
$$
Thus in both cases we get, using,
\begin{eqnarray}
\label{eq1.92}
\int\!\!\!\!\!\!- \ \a^{*\ell} \b^0 F \a^{\ell} \, \vert D \vert^{-1} 
&= &\int\!\!\!\!\!\!- \ \b^0 F \, \vert D \vert^{-1} - \ell 
\int\!\!\!\!\!\!- \ \b_0 F \, \vert D \vert^{-2} \nonumber \\
&= &- \int_{\Hc'}\!\!\!\!\!\!\!\!\!\!- \ \ \vert D \vert^{-1} + \ell 
\int_{\Hc'}\!\!\!\!\!\!\!\!\!\!- \ \ \vert D \vert^{-2}
\end{eqnarray}
the formula,
\be
\label{eq1.93}
\int\!\!\!\!\!\!- \ \mu' [D,\mu] \, \vert D \vert^{-1} = (\ell-k) 
\int_{\Hc'}\!\!\!\!\!\!\!\!\!\!- \ \ \vert D \vert^{-1} + \ell 
(k-\ell)  \int_{\Hc'}\!\!\!\!\!\!\!\!\!\!- \ \ \vert D \vert^{-2} \, .
\ee
One has
$$
\d ([D,\mu]) = \d ((k-\ell) \, \a^{*k} \b^0 F \a^{\ell}) = (k-\ell)^2 
\a^{*k} \b^0 F \a^{\ell}
$$
so that
\be
\label{eq1.94}
\int\!\!\!\!\!\!- \ \mu' \d ([D,\mu]) \, \vert D \vert^{-2} = - 
(k-\ell)^2 \int_{\Hc'}\!\!\!\!\!\!\!\!\!\!- \ \ \vert D \vert^{-2} \, .
\ee
Thus
\begin{eqnarray}
\vp_1 (\mu' , \mu) &= &(\ell-k) \int_{\Hc'}\!\!\!\!\!\!\!\!\!\!- \ \ 
\vert D \vert^{-1} + \left( \ell (k-\ell) + \frac{1}{2} \, (k-\ell)^2 
\right) \int_{\Hc'}\!\!\!\!\!\!\!\!\!\!- \ \ \vert D \vert^{-2} 
\nonumber \\
&= &(\ell-k) \int_{\Hc'}\!\!\!\!\!\!\!\!\!\!- \ \ \vert D \vert^{-1} + 
\frac{1}{2} \, (k^2 - \ell^2) \int_{\Hc'}\!\!\!\!\!\!\!\!\!\!- \ \ 
\vert D \vert^{-2} \, . \nonumber 
\end{eqnarray}
Now one has $\tau_1 
(\mu',\mu)=0$ and, 
$$
b \vp_0 (\mu',\mu) = \vp_0 (\mu'\mu) - \vp_0 (\mu\mu') = \vp_0 
(\a^{*\ell} \b^0 \a^{\ell}) - \vp_0 (\a^{*k} \b^0 \a^k) = \rho (\ell) - 
\rho (k) \, .
$$

 Thus we check that,
\be
\label{eq1.95}
\vp_1 (\mu' , \mu) =  \tau_1 (\mu' , \mu) + b \vp_0 (\mu',\mu) \, .
\ee
Let us now assume that $\mu'$ belongs to the ideal $\Jc$.
We take $\mu' = \a^{*k'} \b^0 \a^{\ell'}$ and $\mu$ to be 
$\a^{*k}$ if $k = \ell' - k' \geq 0$ and $\a^{\ell}$ if $\ell = k' - 
\ell' > 0$. Assume first $k \geq 0$. One has $[D,\mu] = k \a^{*k} F$ 
and $\int\!\!\!\!\!- \ \mu' [D,\mu] \, \vert D \vert^{-1} = k 
\int\!\!\!\!\!- \ \a^{*k'} \b^0 \a^{k'} F \, \vert D \vert^{-1}$. 
Thus using (\ref{eq1.92}) we get,
$$
\int\!\!\!\!\!\!- \ \mu' [D,\mu] \, \vert D \vert^{-1} = k \left( - 
\int_{\Hc'}\!\!\!\!\!\!\!\!\!\!- \ \ \vert D \vert^{-1} + k' 
\int_{\Hc'}\!\!\!\!\!\!\!\!\!\!- \ \ \vert D \vert^{-2} \right) \, .
$$
But $k = \ell' - k'$ so that,
\be
\label{eq1.96}
\int\!\!\!\!\!\!- \ \mu' [D,\mu] \, \vert D \vert^{-1} = (k' - \ell') 
\int_{\Hc'}\!\!\!\!\!\!\!\!\!\!- \ \ \vert D \vert^{-1} + k' (\ell' - 
k') \int_{\Hc'}\!\!\!\!\!\!\!\!\!\!- \ \ \vert D \vert^{-2} \, .
\ee
Also $\d ([D,\mu]) = k^2 \a^{*k} F$ and
$$
\int\!\!\!\!\!\!- \ \mu' \d ([D,\mu] ) \, \vert D \vert^{-2} = k^2 
\int\!\!\!\!\!\!- \ \a^{*k'} \b^0 \a^{k'} F \, \vert D \vert^{-2} = - 
k^2 \int_{\Hc'}\!\!\!\!\!\!\!\!\!\!- \ \ \vert D \vert^{-2} \, .
$$
Thus
$$
\vp_1 (\mu',\mu) = (k'-\ell') \int_{\Hc'}\!\!\!\!\!\!\!\!\!\!- \ \ 
\vert D \vert^{-1} + \left( k' (\ell' - k') + \frac{1}{2} \, k^2 
\right) \int_{\Hc'}\!\!\!\!\!\!\!\!\!\!- \ \ \vert D \vert^{-2} \, .
$$
One has
$$
k' (\ell'-k') + \frac{1}{2} \, k^2 = k' (\ell' - k') + \frac{1}{2} \, 
(\ell'-k')^2 = \frac{1}{2} \, \ell'^2 - \frac{1}{2} \, k'^2 \, ,
$$
so that,
\be
\label{eq1.97}
\vp_1 (\mu' , \mu) = (k'-\ell') \int_{\Hc'}\!\!\!\!\!\!\!\!\!\!- \ \ 
\vert D \vert^{-1} + \frac{1}{2} \, (\ell'^2 - k'^2) 
\int_{\Hc'}\!\!\!\!\!\!\!\!\!\!- \ \ \vert D \vert^{-2} \, .
\ee
Now $\mu'\mu = \a^{*k'} \b^0 \a^{k'}$, $\mu\mu' = \a^{*\ell'} \b^0 
\a^{\ell'}$ so that,
$$
b \vp_0 (\mu',\mu) = \vp_0 (\mu'\mu) - \vp_0 (\mu\mu') = \vp_0 
(\a^{*k'} \b^0 \a^{k'}) - \vp_0 (\a^{*\ell'} \b^0 
\a^{\ell'}) = \rho (k') - 
\rho (\ell') \, .
$$
Thus, since $\tau_1 (\mu',\mu)=0$, we get,
\be
\label{eq1.98}
\vp_1 (\mu',\mu) = \tau_1 (\mu',\mu) + (b \vp_0)(\mu',\mu) \, .
\ee
Next, let us assume that $\ell = k'-\ell' > 0$. Then $\mu = \a^{\ell}$, 
$[D,\mu] = -\ell \a^{\ell} F$, and $\mu' [D,\mu] = -\ell \a^{*k'} \b^0 
\a^{k'} F$ so that by (\ref{eq1.92}),
$$
\int\!\!\!\!\!\!- \ \mu' [D,\mu] \, \vert D \vert^{-1} = - \ell \left( 
- \int_{\Hc'}\!\!\!\!\!\!\!\!\!\!- \ \ \vert D \vert^{-1} + k' 
\int_{\Hc'}\!\!\!\!\!\!\!\!\!\!- \ \ \vert D \vert^{-2} \right) \, ,
$$
\be
\label{eq1.99}
\int\!\!\!\!\!\!- \ \mu' [D,\mu] \, \vert D \vert^{-1} = (k'-\ell') 
\int_{\Hc'}\!\!\!\!\!\!\!\!\!\!- \ \ \vert D \vert^{-1} + k' (\ell' - 
k') \int_{\Hc'}\!\!\!\!\!\!\!\!\!\!- \ \ \vert D \vert^{-2} \, .
\ee
One has $\d ([D,\mu]) = \ell^2 \a^{\ell} F$ and,
\begin{eqnarray}
\int\!\!\!\!\!\!- \ \mu' \d ([D,\mu]) \, \vert D \vert^{-2} &= &\ell^2 
\int\!\!\!\!\!\!- \ \a^{*k'} \b^0 \a^{k'} F \, \vert D \vert^{-1} 
\nonumber \\
&= &-\ell^2 \int_{\Hc'}\!\!\!\!\!\!\!\!\!\!- \ \ \vert D \vert^{-2}  = 
-(k'-\ell')^2 \int_{\Hc'}\!\!\!\!\!\!\!\!\!\!- \ \ \vert D \vert^{-2} \, 
. \nonumber
\end{eqnarray}
Thus as above the coefficient of 
$\int_{\Hc'}\!\!\!\!\!\!\!\!\!\!- \ \ \vert D \vert^{-2}$
in $\vp_1 (\mu' , \mu)$ is  $k' (\ell' - k') + 
\frac{1}{2} \, (k'-\ell')^2 = \frac{1}{2} \, \ell'^2 - \frac{1}{2} \, 
k'^2$ and,
\be \label{eq1.100}
\vp_1 (\mu' , \mu) = (k'-\ell') \int_{\Hc'}\!\!\!\!\!\!\!\!\!\!- \ \ 
\vert D \vert^{-1} + \frac{1}{2} \, (\ell'^2 - k'^2) 
\int_{\Hc'}\!\!\!\!\!\!\!\!\!\!- \ \ \vert D \vert^{-2} \, .
\ee
Thus, as above we get,
\be
\label{eq1.101}
\vp_1 (\mu',\mu) = \tau_1 (\mu',\mu) + b \vp_0 (\mu',\mu) \, .
\ee
Before we proceed, note that (\ref{eq1.89}) which defines
$\vp_0$ is only determined up to the addition of an 
arbitrary constant to $\rho$. As it turns out this
constant will play a role and will be uniquely specified by  equation (\ref{eq1.272}) with the value $\frac{2}{3}$.
Also in order to show that the above computation was largely independent
of the specific numerical values of $\int_{\Hc'}\!\!\!\!\!\!\!\!\!\!- \ \ 
\vert D \vert^{-1}$ and $\int_{\Hc'}\!\!\!\!\!\!\!\!\!\!- \ \ 
\vert D \vert^{-2}$ we did not replace these expressions by their values
which are,
\be
\label{eq1.104}
\int_{\Hc'}\!\!\!\!\!\!\!\!\!\!- \ \ \vert D \vert^{-1} = -1 \, \qquad
\int_{\Hc'}\!\!\!\!\!\!\!\!\!\!- \ \ \vert D \vert^{-2} = 2 \, .
\ee
(To get (\ref{eq1.104}) we use (\ref{eq1.21}) and (\ref{eq1.22}) and 
compute
\begin{eqnarray}
{\rm Trace}_{\Hc'} (e^{-t \vert D \vert}) &= &\left( \sum_{k \in \Zb} 
e^{-t \vert k \vert} \right) \left( \sum_1^{\ify} e^{-t\ell} \right) 
\nonumber \\
&= &\left( 1 + \frac{2}{e^t - 1} \right) \left( \frac{1}{e^t - 1} 
\right) \sim  \frac{2}{t^2} - \frac{1}{t} + \frac{1}{3} - \frac{t}{12} 
+ \cdots) \nonumber 
\end{eqnarray}

\noindent Thus (up to an additive constant), (\ref{eq1.90}) gives,
\be
\label{eq1.105}
\rho (j) = -(j+j^2)
\ee

\noindent We extend the definition of $\vp_0$ to $\Ac$ by 
$\vp_0 (1) = 0$ while, as above, $\vp_0(a)$ 
vanishes if the bidegree of $a$ is $\ne (0,0)$.

\medskip
\subsection{Three dimensional components}

\noindent It follows from (\ref{eq1.101}) that $\psi = \vp_1 - \tau_1 - 
b \vp_0$ vanishes if one of the arguments is in $\Jc$ and thus
$\psi (a_0 , a_1)$ only 
depends on the symbols $\s (a_i)\in C^{\ify}(S^1)$
 where
\be
\label{eq1.51}
0 \longra \Jc \longra \Ac \build\longra_{}^{\s} C^{\ify} (S^1) \longra 
0
\ee
is the natural exact sequence, with $\s (\a) = u$ and $\s (\b) = 0$.

\noindent But the same holds for the component 
$\vp_3$, 
\be
\label{eq1.50}
\vp_3 (a_0 , a_1 , a_2 , a_3) = \int\!\!\!\!\!\!- \ a_0 \, [D,a_1] \, 
[D,a_2] \, [D,a_3] \, \vert D \vert^{-3} \, .
\ee
Indeed if one of the $a_j$ belongs to the two sided ideal $\Jc$   one is dealing with a trace class operator since 
$\vert D \vert^{-3}$ is trace class on the support of $\b$. Thus
$\vp_3 (a_0 , a_1 , a_2 , a_3)$  only depends on the symbols $f_j=\s (a_j)$,
and is given by,
\be
\label{eq1.65}
\vp_3 (a_0 , a_1 , a_2 , a_3) = \frac{-1}{2 \pi i} \int_{S^1} f_0 f'_1 
f'_2 f'_3 {\rm d} \t \, ,
\ee
where $f' = \frac{\partial}{\partial \t}$. 
Since $F = 2P-1 $ introduces a minus sign, we use (\ref{eq1.75})
and can replace $[D,a_j]$ by $- i f'_j$, so that (\ref{eq1.65})
follows from,
\be
\label{eq1.64}
\int\!\!\!\!\!\!- \ \vert D \vert^{-3} = 1 \, .
\ee
Thus to get the complete control of the cocycle $\vp$
it remains only to compute $\psi (\a^k , \a^{*k})$
and $\psi (\a^{*\ell} , \a^{\ell})$.

\noindent Let us compute $\psi (\a^k , \a^{*k})$. One has 
$$
b \vp_0 (\a^k , \a^{*k}) = \vp_0 (\a^k 
\a^{*k}) - \vp_0 (\a^{*k} \a^k)= \, -  \,\vp_0 (\a^{*k} \a^k)
$$
Let $\lb_k = \vp_0 (\a^{*k} \a^k)$, then 
$$
\lb_k - \lb_{k-1} = \vp_0 (\a^{*k-1} (\a^* \a - 1) \, \a^{k-1}) = - 
\vp_0 (\a^{*k-1} \b^0 \a^{k-1}) = - \rho (k-1)
$$
since $\a^* \a - 1 = - \b^0 $. We get,

\be
\label{eq1.103}
\lb_k = \sum_0^{k-1} (j+j^2) = \frac{k^3}{3} - \frac{k}{3} \, ,
\ee
and,
\be
\label{eq.bound}
b \vp_0 (\a^k , \a^{*k}) =\,- \frac{k^3}{3} + \frac{k}{3} 
\ee
Let us compute $\vp_1 (\a^k , \a^{*k})$. With $ \vp_1 = \vp_1^{(0)} - \frac{1}{4} \, \vp_1^{(1)} 
+ \frac{1}{8} \, \vp_1^{(2)}$, one has,
$$
\vp_1^{(0)} (\a^k , \a^{*k}) = \int\!\!\!\!\!\!- \ \a^k [D,\a^{*k}] \, 
\vert D \vert^{-1} = \int\!\!\!\!\!\!- \ k \a^k \a^{*k} F \, 
\vert D \vert^{-1} = k \int\!\!\!\!\!\!- \ F \, \vert D \vert^{-1} 
$$
using (\ref{eq1.75}). Next,
$$
\vp_1^{(1)} (\a^k , \a^{*k}) = \int\!\!\!\!\!\!- \ \a^k \nabla 
([D,\a^{*k}]) \, \vert D \vert^{-3} = k \int\!\!\!\!\!\!- \ \a^k \nabla 
(\a^{*k} F) \, \vert D \vert^{-3} \, .
$$
But $\nabla (x) = \vert D \vert \, \d (x) + \d (x) \, \vert D \vert = 
\d^2 (x) + 2\d (x) \, \vert D \vert$ and $\d^2 (\a^{*k} F) = k^2 
\a^{*k} F$, thus,
$$
\vp_1^{(1)} (\a^k , \a^{*k}) = k^3 \int\!\!\!\!\!\!- \ F \, \vert D 
\vert^{-3} + 2k^2 \int\!\!\!\!\!\!- \ F \, \vert D \vert^{-2}
$$
Finally,
\begin{eqnarray}
\vp_1^{(2)} (\a^k , \a^{*k}) &= &\int\!\!\!\!\!\!- \ \a^k \nabla^2 ([D 
, \a^{*k}]) \, \vert D \vert^{-5} \nonumber \\
&= &4 \int\!\!\!\!\!\!- \ \a^k \d^2 ([D , \a^{*k}]) \, \vert D 
\vert^{-3} = 4 k^3 \int\!\!\!\!\!\!- \ F \, \vert D \vert^{-3} \, . 
\nonumber
\end{eqnarray}
Thus,
\begin{eqnarray}
\label{eq1.108}
\vp_1 (\a^k , \a^{*k}) &= &k \int\!\!\!\!\!\!- \ F \, \vert D 
\vert^{-1} - \frac{1}{2} \, k^2 \int\!\!\!\!\!\!- \ F \, \vert D 
\vert^{-2} \nonumber \\
&+ &\left( - \frac{1}{4} \, k^3 + \frac{1}{2}  \, k^3\right) 
\int\!\!\!\!\!\!- \ F \, \vert D \vert^{-3} \, .
\end{eqnarray}
One has $\int\!\!\!\!\!- \ F \, \vert D \vert^{-3} = -1$ and thus the 
term in $k^3$ is $- \frac{k^3}{4}$. Thus the term in $k^3$ in 
$\psi = \vp_1 - \tau_1 - 
b \vp_0$ 
is, using (\ref{eq.bound}),

\be
\label{eq1.109}
\left( \frac{1}{3} - \frac{1}{4} \right) k^3 = \frac{k^3}{12} \, .
\ee
As we shall see now, this $\frac{1}{12}$ corresponds exactly to
the coefficient $\frac{1}{12}$ in the universal index formula
(theorem 1).

\noindent  Indeed the term in $k^3$ corresponds
to the cochain $\psi_3$ given   in terms of the symbols $f_0 , f_1$ by,
\be
\label{eq1.67}
\psi_3 (a_0 , a_1) = \frac{1}{2\pi i} \int f_0 
f'''_1 {\rm d} \t \, .
\ee
Let us compute $b \psi_3$ where we only involve the symbols. One 
has
$$
f_0 f_1 f'''_2 - f_0 (f_1 f_2)''' + f_2 f_0 f'''_1 = f_0 (-3 f''_1 
f'_2 - 3 f'_1 f''_2) 
$$
thus we get,
\be
\label{eq1.68}
b \psi_3 (a_0 , a_1 , a_2) = \frac{1}{2\pi i} 
\int f_0 (-3 f''_1 f'_2 - 3 f'_1 f''_2) \, {\rm d} \t \, .
\ee
We have,
$$
B_0 \vp_3 (a_0 , a_1 , a_2) = - \frac{1}{2\pi i} \int f'_0 f'_1 f'_2 
{\rm d} \t = \frac{1}{2\pi i} \int f_0 (f''_1 f'_2 + f'_1 f''_2) \, 
{\rm d} \t \, .
$$
This is already cyclic so that,
\be
\label{eq1.69}
B \vp_3 (a_0 , a_1 , a_2) = \frac{3}{2\pi i} \int f_0 (f''_1 f'_2 + 
f'_1 f''_2) \, {\rm d} \t \, .
\ee
and we get,
\be
\label{eqcob}
b \psi_3 \, +  \,B \vp_3= \,0 .
\ee
In fact,
\be
\label{eqbbb}
 \psi_3 \, =  \,B \vp_2  \, \,, \, \,\vp_3= \,b \vp_2 \,,
\ee
where $\vp_2$ is given by,
\be
\label{eqbbc}
\vp_2 (a_0 , a_1 , a_2)=\frac{-1}{2}\frac{1}{2\pi i} \, \int f_0 
f'_1 f''_2 \, {\rm d} \t .
\ee

\noindent Let us now 
compute $\int\!\!\!\!\!- \ F \, \vert D \vert^{-2}$. One has,
\be
\label{eq1.115}
\int\!\!\!\!\!\!- \ \vert D \vert^{-2} = 2 \, .
\ee
and,
\be
\label{eq1.116} \int\!\!\!\!\!\!- \ P \, \vert D \vert^{-2} = 1 \, .
\ee
Thus we get, since $F= 2 P -1$,
\be
\label{eq1.117}
\int\!\!\!\!\!\!- \ F \, \vert D \vert^{-2} = 0 \, .
\ee
and by a similar computation,
\be
\label{eq1.1171}
\int\!\!\!\!\!\!- \ F \, \vert D \vert^{-1} = 1 \, .
\ee
This gives,
\be
\label{eq1.118}
\psi (\a^k , \a^{*k}) = 
\frac{2 k}{3} + \frac{k^3}{12} \, .
\ee
The computation of $\psi (\a^{\ell*} ,\a^{\ell})$ is entirely
similar and gives $\psi (\a^{\ell*} ,\a^{\ell}) = -
\frac{2 \ell}{3} - \frac{\ell^3}{12}$. 
We thus have $\psi=- \frac{2 }{3} \psi_1 + \frac{1}{12} \psi_3$ where,

\be
\label{eq1.167}
\psi_1 (a_0 , a_1) = \frac{1}{2\pi i} \int f_0 
f'_1 {\rm d} \t \, .
\ee
It just remains to see why adding a constant to $\rho$ allows
to eliminate $- \frac{2 }{3} \psi_1 $ from $ \psi$. This
follows from (\ref{eq1.103}) and (\ref{eq.bound}) i. e.
\be
\label{eq1co}
b \vp_0 (\a^k , \a^{*k}) =\, \sum_0^{k-1} \rho(j) \, ,
\ee
Thus adding $\frac{2 }{3}$ to $\rho$ gives $\psi=\frac{1}{12} \psi_3$
and ends the proof of theorem 3. $\Box$

\noindent We shall now understand the conceptual meaning of the 
above concrete computation.

\bigskip
\section{The $\eta$-Cochain.}

In this section we shall give two general formulas. The first
will provide the conceptual explanation of theorem 3, and of the 
cochain ($\vp_0,\vp_2$) which appears there. The second will
prepare for the computation of the local index formula
in the general case $q \in ]0,1[$.

\noindent 
The explanation of theorem 3 and of the cochains, 
\begin{eqnarray}
\label{phizero}
&&\vp_0 (\a^{*j} f(\b)
\,\a^j) =  (\frac{2}{3} - j - j^2)\, \, \frac{1}{2\pi}  \int_{S^1} f \, \,{\rm d} \t ,\\
&& \vp_2(a_0,a_1,a_2)=\,\frac{-1}{24}\frac{1}{2\pi i} \, \int f_0 
f'_1 f''_2 \, {\rm d} \t, \,\,\,f_j= \s(a_j),
\end{eqnarray}
is given by the following,

\bigskip

\noindent {\bf Proposition 2.} {\it
Let $(\Ac ,\Hc ,D)$ be a spectral triple with discrete simple dimension
spectrum not containing $0$ and upper bounded by $ 3$. Assume
that $[ F, a]$ is trace class for all $a \in \Ac$.
Let $\tau_1(a_0,a_1)= {\rm Trace}(a_0 [ F, a_1])$. 

\noindent Then the local Chern Character ($\vp_1,\vp_3$) of $(\Ac ,\Hc ,D)$
is equal to  $\tau_1 + (b+B)\vp$ where ($\vp_0,\vp_2$) is the cochain given by,
}
\begin{eqnarray}
\label{prop2}
&&\vp_0 (a) =  {\rm Trace}(F \,a\, \vert D \vert^{-s})_{s=0},\nonumber \\
&& \vp_2(a_0,a_1,a_2)=\,\frac{1}{24}\int\!\!\!\!\!\!- \ a_0  \,\d(a_1) \,\d^2( a_2) \,F  \, \vert 
D \vert^{-3} \,.\nonumber
\end{eqnarray}

\bigskip
\noindent Note that $\vp_0$ makes sense by the absence of pole at $s=0$,
i.e. the hypothesis $0 \notin {\rm Dimension \,\,Spectrum}$.
Its value for $a=1$ coincides with the
classical $\eta$-invariant (\cite{APS1},\cite{APS2}) and justifies the terminology of $\eta$-cochain to qualify
the cochain ($\vp_0,\vp_2$). 

\noindent The proof of the proposition is a simple calculation based on the expansion
(\cite{CM}) 
\begin{eqnarray}
\label{newton}
&&\vert D \vert^{-s} a \sim a \,\,\vert D \vert^{-s} -s \ \d (a)\vert D \vert^{-s-1} + {(-s)(-s-1) \over 2!} \ \d^2 (a) \vert D \vert^{-s-2}\nonumber\\
&&+ \,\,{(-s)(-s-1)(-s-3) \over 3!} \ \d^3 (a) \vert D \vert^{-s-3} + s\,\, \,o(\vert D \vert^{-s-3}),
\end{eqnarray}
which allows to express $b \vp_0$ in terms of residues. More specifically
one gets,

\begin{eqnarray}
\label{bphi}
&&b\vp_0 (a_0,a_1) = -\tau_1(a_0,a_1) +  \,\int\!\!\!\!\!\!- \ a_0  \,\d(a_1) \,F  \, \vert 
D \vert^{-1} \,\nonumber\\
&&-\,\frac{1}{2} \,\int\!\!\!\!\!\!- \ a_0  \,\d^2(a_1) \,F  \, \vert 
D \vert^{-2} +\,\frac{1}{3} \,\int\!\!\!\!\!\!- \ a_0  \,\d^3(a_1) \,F  \, \vert 
D \vert^{-3} \,
\end{eqnarray}
using the hypothesis $[ F, a]$ trace class for all $a \in \Ac$.
This hypothesis also shows that,
\begin{eqnarray}
\label{phi1}
&&\vp_1 (a_0,a_1) =   \,\int\!\!\!\!\!\!- \ a_0  \,\d(a_1) \,F  \, \vert 
D \vert^{-1} \,\nonumber\\
&&-\,\frac{1}{2} \,\int\!\!\!\!\!\!- \ a_0  \,\d^2(a_1) \,F  \, \vert 
D \vert^{-2} +\,\frac{1}{4} \,\int\!\!\!\!\!\!- \ a_0  \,\d^3(a_1) \,F  \, \vert 
D \vert^{-3} \,
\end{eqnarray}
Comparing (\ref{bphi}) with (\ref{phi1}) gives the required $\frac{1}{12}$
and allows to check that $\vp_{odd}=\tau_1 + (b+B)\vp_{ev}$.

\noindent Let us compute  $\vp_{ev}$ in the above example. One has, as in (\ref{eq2.47}),
\be
\label{example1}
\vp_0 (\a^{*k}e
\,\a^k)=({\rm Trace} \, (F\,e \, (\vert D \vert +k)^{-s}))_{s=0} 
\ee
Using $F= 2 P -1$ this gives,
\be
\label{example2}
\vp_0 (\a^{*k} e
\,\a^k)= 2(\sum_0^{\ify} \frac{1}{(n+k)^s} )_{s=0} 
-(\sum_0^{\ify} \frac{2n+1}{(n+k)^s})_{s=0} 
\ee
One has,
$$
(\sum_0^{\ify} \frac{1}{(n+k)^s} )_{s=0}=\Zc(0)-(k-1).
$$
Also,
\begin{eqnarray}
\label{zeta}
&&\sum_0^{\ify} \frac{2n+1}{(n+k)^s})_{s=0} = \, \,
2\Zc(-1)+(1-2 k) \,\Zc(0)-\sum_1^{k-1}(2 \ell -2 k +1)\nonumber \\
&&= \, \, 2 \,\Zc(-1)+(1-2 k) \,\Zc(0) + (k-1)^2.
\end{eqnarray}
Thus we get,
\be
\label{example3}
\vp_0 (\a^{*k}e
\,\a^k)=2 (\Zc(0)-(k-1))-( 2 \,\Zc(-1)+(1-2 k) \,\Zc(0) + (k-1)^2)
\ee
which using the values, 
\be
\label{zeta1}
\Zc(0)=-\frac{1}{2} \, ,\,\,\,\,\,\Zc(-1)=-\frac{1}{12}
\ee
gives the desired result,
\be
\label{example3}
\vp_0 (\a^{*k}e
\,\a^k)=\frac{2}{3} - k - k^2.
\ee
The only other non-trivial value of $\vp_0$
is $\eta= \vp_0(1)$, and the computation gives $\eta= \frac{1}{2}$.
Finally the equality,
\be
\int\!\!\!\!\!\!- \ a_0  \,\d(a_1) \,\d^2( a_2) \,F  \, \vert 
D \vert^{-3} =\,\,\frac{1}{2\pi i} \, \int f_0 
f'_1 f''_2 \, {\rm d} \t, \,\,\,f_j= \s(a_j),
\ee
and the coincidence of the functional $\tau_1$
of theorem 3 with $ {\rm Trace}(a_0 [ F, a_1])$ give a perfect
account of theorem 3.

\noindent In order to lighten the general computation, for $q \in ]0,1[$,
we shall state a small variant of proposition 2, proved in a similar way.
Given a spectral triple $(\Ac ,\Hc ,D)$ let us define the metric dimension Dm$(P)$ of 
a projection $P$ commuting with $D$ as the lower 
bound of all $d \in \Rb$ such that $P (D+i)^{-1}$ is in the Schatten class
$L^d$. We then have as above,

\bigskip

\noindent {\bf Proposition 3.} {\it
Let $(\Ac ,\Hc ,D)$ be a spectral triple with discrete dimension
spectrum not containing $0$. Assume that} Dm$(\Hc)\leq 3$, {\it and} Dm$(P)$
 $\leq 2$, $P=(1+F)/2$,
{ \it and that $[ F, a]$ is trace class for all $a \in \Ac$. 
\noindent Then the local Chern Character ($\vp_1,\vp_3$) of $(\Ac ,\Hc ,D)$
is equal to  $\psi_1 -(b+B)\vp$ where $\psi_1$ is the cyclic cocycle,} 
$$
\psi_1(a_0,a_1)=2 \int\!\!\!\!\!\!- \ a_0  \,\d(a_1) \, P\, \vert 
D \vert^{-1} \,- \int\!\!\!\!\!\!- \ a_0  \,\d^2(a_1) \, P\, \vert 
D \vert^{-1} \,
$$ 
 {\it and ($\vp_0,\vp_2$) is the cochain given by,
}
\begin{eqnarray}
\label{prop3}
&&\vp_0 (a) =  {\rm Trace}(\,a\, \vert D \vert^{-s})_{s=0},\nonumber \\
&& \vp_2(a_0,a_1,a_2)=\,\frac{1}{24}\int\!\!\!\!\!\!- \ a_0  \,\d(a_1) \,\d^2( a_2) \,  \, \vert 
D \vert^{-3} \,.\nonumber
\end{eqnarray}

\medskip \noindent Combining Propositions 2 and 3 one obtains
under the hypothesis of Proposition 3 the equality,
\be
\label{eta2}
\psi_1 - \tau_1= \, b (\psi_0)
\ee
where the cochain $\psi_0$ is given by,
\be
\label{eta3}
\psi_0(a)= \, 2\; {\rm Trace}(\,a\, P \,
\vert
 D \vert^{-s})_{s=0}
\ee

\bigskip
\section{Pseudo-differential calculus and the cosphere bundle on ${\rm SU}_q (2)$
, $q \in \, ] 0,1 [$}

In this section we shall construct the pseudo-differential calculus on 
${\rm SU}_q (2)$ following the general theory of \cite{CM}. We shall 
determine the algebra of complete symbols by computing the quotient by 
smoothing operators. This will give the cosphere bundle $S_q^*$ 
of ${\rm SU}_q (2)$ and the analogue of the 
geodesic flow will yield a one-parameter group of
automorphisms $\g_t$ of $C^{\ify} (S_q^*)$. We shall also
construct the restriction morphism $r$ to the 
product of two 2-disks,
 
\be
\label{eq4.1}
r : C^{\ify} (S_q^*) \ra C^{\ify} (D_{q_+}^2 \ts D_{q_-}^2)
\ee

\smallskip

\noindent Our goal is to prepare for the computation in the next section 
of the dimension spectrum and of residues. Let us recall from \cite{CM} 
that given a spectral triple $(\Ac , \Hc , D)$ we say that an operator $P$ 
in $\Hc$ is of order $\a$ when,
\be
\label{eq3.1}
\vert D \vert^{-\a} \, P \in \bigcap_{n=1}^{\ify} {\rm Dom} \, \d^n
\ee
where $\d$ is the unbounded derivation given by,
\be
\label{eq3.2}
\d (T) = \vert D \vert \, T - T \, \vert D \vert \, .
\ee
Thus $OP^0 = \build\bigcap_{n=1}^{\ify} {\rm Dom} \, \d^n$ is the algebra 
of operators of order $0$ and $OP^{-\ify} = \build\bigcap_{k > 0}^{} 
OP^{-k}$ is a two sided ideal in $OP^0$.

\smallskip

\noindent We let $(\Ac , \Hc , D)$ be the spectral triple of \cite{CP} and we 
first determine the algebra $\Bc$ generated by the $\d^k (a)$, $a \in 
\Ac$.

\noindent  Recall that $D$ is the diagonal operator in $\Hc$ given by,
\be
\label{eq2.18}
D (e_{ij}^{(n)}) = (2 \, \d_0 (n-i)-1) \, 2n \,\, e_{ij}^{(n)}
\ee
where $\d_0 (k) = 0$ if $k \ne 0$ and $\d_0 (0) = 1$.  
\smallskip

\noindent By construction, the generators $\a , \b$ of $\Ac$ are of the 
form,
\be
\label{eq3.3}
\a = \a_+ + \a_- \, , \qquad \b = \b_+ + \b_-
\ee
where,
\be
\label{eq3.4}
\d (\a_{\pm}) = \pm \a_{\pm} \, , \qquad \d (\b_{\pm}) = \pm \b_{\pm} \, .
\ee
The explicit form of $\a_{\pm}$, $\b_{\pm}$ is, using $\frac{n}{2}$ 
instead of $n$ for the notation of the $\frac{1}{2}$ integer,
\be
\label{eq3.5}
\a_{\pm} (e_{(i,j)}^{(n/2)}) = a_{\pm} \left( n/2 , i , j \right) 
\, e_{\left( i - \frac{1}{2} , j - \frac{1}{2} \right)}^{\left( \frac{n 
\pm 1}{2} \right)}
\ee
\be
\label{eq3.6}
\b_{\pm} (e_{(i,j)}^{(n/2)}) = b_{\pm} \left( n/2 , i , j \right) 
\, e_{\left( i + \frac{1}{2} , j - \frac{1}{2} \right)}^{\left( \frac{n 
\pm 1}{2} \right)}
\ee
where $a_{\pm}$, $b_{\pm}$ are as in (\ref{eq25}) and (\ref{eq26}) above.

\smallskip

\noindent Thus the algebra $\Bc$ is generated by the operators $\a_{\pm}$, 
$\b_{\pm}$ and their adjoints.

\smallskip

\noindent We shall now see that, modulo the smoothing operators, we can 
strip the complicated formulas for the coefficients $a_{\pm}$, $b_{\pm}$ 
and replace them by extremely simple ones. Since we are computing {\it 
local formulas} we are indeed entitled to mod out by smoothing operators 
and this is exactly where great simplifications do occur.

\smallskip

\noindent Let us first relabel the indices $i,j$ using,
\be
\label{eq3.7}
x = \frac{n}{2} + i \, , \qquad y = \frac{n}{2} + j \, .
\ee
By construction $x$ and $y$ are {\it integers} which vary exactly in $\{0 
, 1 , \ldots , n \}$.

\smallskip

\noindent Working modulo $OP^{-\ify}$ means that we can neglect in the 
formulas for $a_{\pm}$, $b_{\pm}$ any modification by a sequence of rapid 
decay in the set:
\be
\label{eq3.8}
\Lb = \{ (n,x,y) \, ; \ n \in \Nb \, , \ x,y \in \{0 , \ldots , n \}\} \, 
.
\ee
Thus first, we can  get rid of the denominators, since both 
$(1-q^{2n})^{-1/2}$ or $(1-q^{(2n+2)})^{-1/2}$ are equivalent to $1$
and the numerators are bounded.

\smallskip

\noindent Next, when we rewrite the numerators in terms of the variables 
$n,x,y$ we get, say for $a_+$, the simplified form,
\be
\label{eq3.9}
a'_+ (n,x,y) = q^{1+x+y} (1-q^{2+2(n-x)})^{1/2} (1-q^{2+2(n-y)})^{1/2} \, 
.
\ee
Modulo sequences of rapid decay one has,
$$
q^x (1-q^{2+2(n-x)})^{1/2} \sim q^x \, ,
$$
as one sees from the inequality $(1-(1-u)^{1/2}) \leq u$ valid for $u \in 
[0,1]$, and the fact that
$$
q^x \, q^{2(n-x)} \leq q^x \, q^{(n-x)} = q^n \, .
$$
Thus we see that modulo sequences of rapid decay we can replace $a'_+$ by,
\be
\label{eq3.10}
a''_+ (n,x,y) = q^{1+x+y} \, .
\ee
To simplify formulas let us relabel the basis as,
\be
\label{eq3.11}
f_{x,y}^{(n)} = e_{\left( x- n/2 ,y - n/2 
\right)}^{(n/2)} \, ,
\ee
then the following operator agrees with $\a_+$ modulo $OP^{-\ify}$,
\be
\label{eq3.12}
\a'_+ (f_{x,y}^{(n)}) = q^{1+x+y} \, f_{x,y}^{(n+1)} \, .
\ee
For $\a_-$ one has, as above,
\be
\label{eq3.13}
a'_- (n,x,y) = (1 - q^{2x})^{1/2} \, (1-q^{2y})^{1/2}
\ee
and the corresponding operator $\a'_-$ is,
\be
\label{eq3.14}
\a'_- (f_{x,y}^{(n)}) = (1-q^{2x})^{1/2} (1-q^{2y})^{1/2} \, 
f_{x-1,y-1}^{(n-1)} \, .
\ee
Note that $\a'_-$ makes sense for $x=0$, $y=0$. For $\b_+$ one gets,
\be
\label{eq3.15}
b'_+ (n,x,y) = -q^y (1-q^{2+2(n-y)})^{1/2} (1-q^{2+2x})^{1/2}
\ee
and as above we can replace it by,
\be
\label{eq3.16}
b''_+ (n,x,y) = -q^y (1-q^{2+2x})^{1/2}
\ee
which gives,
\be
\label{eq3.17}
\b'_+ (f_{x,y}^{(n)}) = -q^y (1-q^{2+2x})^{1/2} \, f_{x+1,y}^{(n+1)} \, .
\ee
In a similar way one gets,
\be
\label{eq3.18}
\b'_- (f_{x,y}^{(n)}) = q^x (1-q^{2y})^{1/2} \, f_{x,y-1}^{(n-1)}
\ee
which makes sense even for $y=0$.

\smallskip

\noindent It is conspicuous in the above formulas that the new and much 
simpler coefficients {\it no longer depend upon the variable} $n$.

\smallskip

\noindent To understand these formulas we introduce the
following representations 
$\pi_{\pm}$ of $\Ac = C^{\ify} ({\rm SU}_q (2))$ \footnote{see the appendix
for the notation}. In both cases the Hilbert spaces are $\Hc_{\pm} = \ell^2 
(\Nb)$ with basis $(\ve_x)_{x \in \Nb}$ and the representations are given 
by,
\be
\label{eq3.19}
\pi_{\pm} (\a) \, \ve_x = (1-q^{2x})^{1/2} \, \ve_{x-1} \qquad \fl x \in 
\Nb
\ee
\be
\label{eq3.20}
\pi_{\pm} (\b) \, \ve_x = \pm \,\,q^x \, \ve_x \qquad \fl x \in \Nb \, .
\ee
With these notations, and if we ignore the $n$-dependence in the above 
formulas we have the correspondence,
\begin{eqnarray}
\label{eq3.21}
\a'_+ \cong -q \, \b^* \ot \b \\
\a'_- \cong \a \ot \a \nonumber \\
\b'_+ \cong \a^* \ot \b \nonumber \\
\b'_- \cong \b \ot \a \, , \nonumber
\end{eqnarray}
through the representation $\pi = \pi_+ \ot \pi_-$.
Now recall that $\Ac$ is a Hopf algebra, with coproduct corresponding to 
matrix tensor multiplication for the following $2 \ts 2$ matrix,
\be
\label{eq3.22}
U = \left[ \matrix{\a &-q \b^* \cr \b &\a^* \cr} \right]
\ee
which gives,
\be
\label{eq3.23}
\D \a = \a \ot \a - q \b^* \ot \b
\ee
$$
\D \b = \b \ot \a + \a^* \ot \b \, .
$$
This shows of course that $\a' = \a'_+ + \a'_-$ and $\b' = \b'_+ + \b'_-$ 
provide a representation of $\Ac$ which is the tensor product
in the sense of Hopf algebras of the representations $\pi_+ $ and $ \pi_-$
of $\Ac$. However to 
really understand the algebra $\Bc$ modulo $OP^{-\ify}$ an its action in 
$\Hc$ we need to keep track of  the shift of $n$ in the formulas
for $\a'_{\pm}$ and $\b'_{\pm}$.

\smallskip

\noindent One can encode these shifts using the $\Zb$-grading of 
$\Bc$ coming from the one parameter group of automorphisms $\g (t)$ which 
plays the role of the geodesic flow,
\be
\label{eq3.24}
\g (t) (P) = e^{it \vert D \vert} \, P \, e^{-it \vert D \vert} \, .
\ee
For the corresponding $\Zb$-grading one has,
\be
\label{eq3.25}
\deg (\a_{\pm}) = \pm 1 \, , \qquad \deg (\b_{\pm}) = \pm 1 \, ,
\ee
which are the correct powers of the shifts of $n$ in the above 
formulas
for $\a'_{\pm},\b'_{\pm}$.
To $\g$ we associate the algebra morphism,
\be
\label{eq3.26}
\g : \Bc \ra \Bc \ot C^{\ify} (S^1) = C^{\ify} (S^1 , \Bc)
\ee
given by $\g (b)(t) = \g_t (b)$, $ \fl t \in S^1$.

\noindent Finally, note that the representations $\pi_{\pm}$ are not faithful
on $C^{\ify} (SU_q(2))$ since the spectrum of $\b$ is real and positive
in $\pi_+$ and real negative for $\pi_-$. We let $C^{\ify} (D_{q\pm}^2)$
be the corresponding quotient algebras and $r_{\pm}$ the restriction
morphisms.
\medskip

\noindent {\bf Proposition 4.} {\it The following equalities define an 
algebra homomorphism $\rho$ from $\Bc$ to
$$
C^{\ify} (D_{q+}^2) \ot C^{\ify} (D_{q-}^2) \ot C^{\ify} (S^1) \, ,
$$
$$
\rho (\a_+) = -q \b^* \ot \b \ot u \, , \quad \rho (\a_-) = \a \ot \a \ot 
u^* \, ,
$$
$$
\rho (\b_+) = \a^* \ot \b \ot u \, , \quad \rho (\b_-) = \b \ot \a \ot u^* 
\, ,
$$
where we omitted $r_+ \ot r_-$.}

\medskip

\noindent {\bf Proof.} Using (\ref{eq3.26}) it is enough to show that the 
formulas,
$$
\rho_1 (\a_+) = -q \, \pi_+ (\b^*) \ot \pi_- (\b) \, , \quad \rho_1 (\a_-) 
= \pi_+ (\a) \ot \pi_- (\a)
$$
$$
\rho_1 (\b_+) = \pi_+ (\a^*) \ot \pi_- (\b) \, , \quad \rho_1 (\b_-) = 
\pi_+ (\b) \ot \pi_- (\a)
$$
define a representation of $\Bc$.

\smallskip

\noindent But this representation is weakly contained in the natural 
representation of $\Bc$ in $\Hc$. To obtain $\rho_1$ from the latter, one 
just considers vectors $\ve_{x,y}^N$ in $\Hc$, of the form,
\be
\label{eq3.27}
\ve_{x,y}^N = \sum h_{(n)}^N \, f_{x,y}^{(n)}
\ee
where $h^N \in \ell^2 (\Nb)$ corresponds to the amenability of the group 
$\Zb$, i.e. to the weak containement of the trivial representation of 
$\Zb$ by the regular one. Thus $h^N$ depends on a large integer $N$ and is 
$1/\sqrt N$ for $0 \leq n < N$ and $0$ for $n \geq N$.

\smallskip

\noindent The almost invariance of $h^N$ under translation of $n$ shows 
that the $n$-dependence of the formulas (\ref{eq3.15})--(\ref{eq3.18}) 
disappears when $N \ra \ify$ and that $\rho_1$ is a representation of $\Bc$. 
Finally $\rho $ is its amplification using (\ref{eq3.26}) $\Box$

\bigskip

\noindent {\bf Definition 1.} {\it Let $C^{\ify} (S_{q}^*)$ be the range 
of $\rho$ in $C^{\ify} (D_{q+}^2 \ts D_{q-}^2 \ts S^1)$.
}
\bigskip

\noindent  By construction  $C^{\ify} (S_{q}^*)$ is topologically
generated by $\rho (\a_{\pm})$, $\rho (\b_{\pm})$. The NC-space
$S_{q}^*$ plays the
role of the  {\it cosphere bundle}. 
 The algebra $C^{\ify} (S_{q}^*)$ is strictly contained in $C^{\ify} (D_{q+}^2 \ts D_{q-}^2 \ts S^1)$ since its image under $\s \ot \s \ot {\rm Id}$ is  the subalgebra of $C^{\ify}( S^1 \ts S^1 \ts S^1)$ generated by $u \ot u \ot u^*$.
Let $ \nu_t$ be the $S^1$-action on $S_{q}^*$ given by 
the restriction of the derivation
$1 \ot 1 \ot {\partial_u}$ where $\partial_u(u)=u$. By construction, \be
\label{eq3.29}
\rho( \g (t) (P))= \nu_t(\rho (P)
\ee
so that $\nu_t$ is the analogue of the action of the 
geodesic flow on the cosphere bundle. We let,
\be
\label{eq4.1}
r : C^{\ify} (S_q^*) \ra C^{\ify} (D_{q_+}^2 \ts D_{q_-}^2)
\ee
be the natural restriction morphism.

\noindent  Viewing $\rho$ as the total symbol map we shall now define a 
natural lifting from symbols to operators. This will only be relevant on 
the range of $\rho$ but to define it we start from
the representation $\pi=\pi_+ \ot \pi_- \ot s$ of $ C^{\ify} (D_{q_+}^2) \ot C^{\ify} 
(D_{q_-}^2) \ot C^{\ify} (S^1)$ in $\ell^2 (\Nb) \ot \ell^2 (\Nb) \ot \ell^2 (\Zb)$ where $s(u)$ is the shift $S$ in $\ell^2 (\Zb)$.
We let  $Q$ be the
 orthogonal projection on the subset $\Lb$ of the basis $f_{x,y}^{(n)}$ 
determined by $n \geq { \rm sup}(x,y)$ and identify the range of $Q$
with the Hilbert space $\Hc$.
By definition the lifting $\lb$ is the compression,
\be
\label{eq3.35}
\lb (g) = Q \, \pi(g)  \, Q
\ee
For $g $ 
 of the form $\mu \ot u^n$, one has,
\be
\label{eq3.33}
\lb (g) \, f_{x,y}^{(\ell)} = \sum \mu_{(x,y)}^{(x',y')} \, 
f_{x',y'}^{(\ell+n)}
\ee
where $\mu_{(x,y)}^{(x',y')}$ are the matrix elements for the action of 
$\mu$ in $\ell^2 (\Nb) \ot \ell^2 (\Nb)$,
\be
\label{eq3.34}
\mu \, \ve_{x,y} = \sum \mu_{(x,y)}^{(x',y')} \, \ve_{x',y'} \, .
\ee
It may happen in formula (\ref{eq3.33}) that the indices in 
$f_{x',y'}^{(\ell+n)}$ do not make sense, i.e. that $f_{x',y'}^{(\ell+n)}$
does not belong to $\Lb$. In that case the corresponding 
term is $0$. We have now restaured the shift of $n$ in the formulas
for $\a'_{\pm}$ and $\b'_{\pm}$ and get,
\medskip

\noindent {\bf Lemma 1.} {\it For any $b \in \Bc$ one has,
$$
b - \lb (\rho (b)) \in OP^{-\ify} \, .
$$
}
\noindent We refer to the appendix for the implications of this lemma. We give there another
general lemma proving the stability under 
holomorphic functional calculus for the natural smooth algebras involved in our discussion.

\bigskip
\section{Dimension Spectrum and Residues for ${\rm SU}_q (2)$,  $q \in 
\, ]0,1[$}

Let as above $S_q^*$ be the cosphere bundle of ${\rm SU}(2)_q$, $\g_t$ its 
geodesic flow and,
\be
\label{eq4.1}
r : C^{\ify} (S_q^*) \ra C^{\ify} (D_{q_+}^2 \ts D_{q_-}^2)
\ee
be the natural restriction morphism.

\smallskip

\noindent For $C^{\ify} (D_{q_{\pm}}^2)$ we have an exact sequence of the 
form,
\be
\label{eq4.2}
0 \longra \Sc \longra C^{\ify} (D_q^2) \build\longra_{}^{\s} C^{\ify} (S^1) 
\longra 0
\ee
where the ideal $\Sc$ is isomorphic to the algebra of matrices of rapid 
decay. Using the representations $\pi_{\pm}$ of $C^{\ify}(D_{q_{\pm}}^2)$ in $\ell^2 
(\Nb)$ with basis $(\ve_x)$, $x \in \Nb$, we define two linear functionals 
$\tau_0$ and $\tau_1$ by,
\be
\label{eq4.3}
\tau_1 (a) = \frac{1}{2 \pi} \int_0^{2\pi} \s (a) \, {\rm d} \t \qquad \fl 
a \in C^{\ify} (D_q^2)
\ee
and
\be
\label{eq4.4}
\tau_0 (a) = \lim_{N \ra \ify} {\rm Trace}_N (\pi(a)) - \tau_1 (a) \, N
\ee
where,
\be
\label{eq4.5}
{\rm Trace}_N (a) = \sum_0^N \langle a \, \ve_x , \ve_x \rangle \, .
\ee
$($where we omitted $\pm$ in the above formulas$)$.
For $a \in \Sc$ one has $\s (a) = 0$ and $\tau_1 (a) = 0$, $\tau_0 (a) = 
{\rm Trace} \, (a)$. In general both $\tau_0$ and $\tau_1$ are invariant 
under the one parameter group generated by $\partial_{\a}$ and on the fixed 
points of this group, one has,
\be
\label{eq4.6}
\tau_0 (a) = {\rm Trace} \, (\pi(a)- \tau_1 (a) \, 1) + \tau_1 (a) \, .
\ee
For all $a \in \Ac$ one has, (for all $k > 0$),
\be
\label{eq4.7}
{\rm Trace}_N (\pi(a)) = \tau_1 (a) \, N + \tau_0 (a) + 0(N^{-k}) \, .
\ee
We shall now prove a general formula 
computing residues
 of pseudo-differential operators in terms of their symbols,

\bigskip

\noindent {\bf Theorem 4.} 

\begin{enumerate}
\item {\it The dimension spectrum of ${\rm SU}_q (2)$ 
is $\{ 1,2,3 \}$.}

\smallskip

\item {\it Let $b \in \Bc$, $\rho (b) \in C^{\ify} (S_q^*)$ its 
symbol. Then let $\rho (b)^0$ be the component of degree $0$ for the 
geodesic flow $\g_t$. One has,
\[
\int\!\!\!\!\!\!- \, b \, \vert D \vert^{-3} = (\tau_1 \ot \tau_1)(r\rho 
(b)^0)
 \]
\[
\int\!\!\!\!\!\!- \, b \, \vert D \vert^{-2} = (\tau_1 \ot \tau_0 + \tau_0 
\ot \tau_1)(r\rho (b)^0)
 \]
\[
\int\!\!\!\!\!\!- \, b \, \vert D \vert^{-1} = (\tau_0 \ot \tau_0)(r\rho 
(b)^0)\, .
 \]
}\end{enumerate}

\medskip

\noindent {\bf Proof.} By lemma the operator $b-\lb \rho (b)$ belongs to $ 
OP^{-\ify}$, thus 

 \noindent  $\Zc_b (s) - {\rm Trace} \, (\lb \rho (b)$ $\vert D 
\vert^{-s})$ is a holomorphic function of $s \in \Cb$.

 \noindent  One has  ${\rm Trace} 
\, (\lb \rho (b) \, \vert D \vert^{-s}) = {\rm Trace} \, (\lb \rho (b)^0 
\vert D \vert^{-s})$ and with $\rho (b)^0 = T$,
\be
\label{eq4.8}
{\rm Trace} \, (\lb (T) \, \vert D \vert^{-s}) = \sum_{n=0}^{\ify} \left( 
\sum_{x=0}^n \ \sum_{y=0}^n \langle \pi (T) \, \ve_{x,y} , \ve_{x,y} 
\rangle \right) n^{-s} \, .
\ee
Thus by (\ref{eq4.7}) we get, modulo holomorphic functions of $s \in \Cb$,
\begin{eqnarray}
\label{eq4.9}
{\rm Trace} \, (\lb (T) \, \vert D \vert^{-s}) &\cong &(\tau_1 \ot 
\tau_1)(T) \, \Zc (s-2) \\
&+ &(\tau_1 \ot \tau_0 + \tau_0 \ot \tau_1)(T) \, \Zc (s-1) + (\tau_0 \ot 
\tau_0)(T) \, \Zc (s) \, . \nonumber
\end{eqnarray}
This shows that $\Zc_b (s)$ extends to a meromorphic function of $s$ with 
simple poles at $s \in \{ 1,2,3 \}$ and gives the above values for the 
residues.

\smallskip

\noindent To show 1) we still need to adjoin $F = {\rm Sign} \, D$ to the 
algebra $\Bc$, but by \cite{CP} one has,
\be
\label{eq4.10}
[F,a] \in OP^{-\ify} \qquad \fl a \in \Bc
\ee
so that the only elements which were not handled above are those of the 
form,
\be
\label{eq4.11}
bP \, , \quad b \in \Bc \, , \quad P = \frac{1+F}{2} \, .
\ee
Thus with the above notation we still need to analyse,
\be
\label{eq4.12}
{\rm Trace} \, (\lb (T) \, P \, \vert D \vert^{-s}) \, .
\ee
Since $P$ corresponds to the subset of the basis $f_{x,y}^{(n)}$ given by 
$\{ x=n \}$ in the above notations, the trace (\ref{eq4.12}) can be 
expressed as,
\be
\label{eq4.13}
{\rm Trace} \, (\lb (T) \, P \, \vert D \vert^{-s}) = \sum_{n=0}^{\ify} \ 
\sum_{y=0}^n \langle \pi (T) \, \ve_{n,y} , \ve_{n,y} \rangle \, n^{-s}
\ee
and the structure of the representation $\pi_+$ shows that the r.h.s. 
gives, modulo holomorphic function of $s \in \Cb$,
\be
\label{eq4.14}
{\rm Trace} \, (\lb (T) \, P \, \vert D \vert^{-s}) \cong (\tau_1 \ot 
\tau_1)(T) \, \Zc (s-1) + (\tau_1 \ot \tau_0) (T) \, \Zc (s) \, .
\ee
This shows (\ref{eq4.1}) and also gives the two formulas,
\be
\label{eq4.15}
\int\!\!\!\!\!\!- \, b \, P \, \vert D \vert^{-2} = (\tau_1 \ot 
\tau_0)(\rho (b)^0)
\ee
and,
\be
\label{eq4.16}
\int\!\!\!\!\!\!- \, b \, P \, \vert D \vert^{-1} = (\tau_0 \ot \tau_0) \, 
\rho (b)^0
\ee
which we shall now exploit to do the computation of the local index formula 
for ${\rm SU}_q (2)$.

\section{The local index formula for ${\rm SU}_q (2)$, $q \in \, ] 
0,1[$}

The local index formula for the spectral triple of ${\rm SU}_q (2)$
uniquely determines a cyclic 1-cocycle and hence by (\cite{Coober2})
a corresponding one dimensional cycle. We shall first
describe independently the obtained cycle since the 
NC-differential calculus it exhibits is of independent interest.

\noindent Let $\Ac = C^{\ify} ({\rm SU}_q (2))$ and $\partial$ the derivation,
\be
\label{eq4.17}
\partial = \partial_{\b} - \partial_{\a} \, .
\ee
We extend the functional $\tau_0$ of (\ref{eq4.4}) to $\Ac$ by,
\be
\label{eq4.18}
\tau (a) = \tau_0 \,( r_- (a^{(0)})) \qquad \fl a \in \Ac
\ee
where $a^{(0)}$ is the component of degree $0$ for $\partial$. By 
construction $\tau$ is $\partial$-invariant but fails to be a trace. It is 
the average of the transformed of $\tau_0 \circ r_-$ by the automorphism 
$\nu_t \in {\rm Aut} (\Ac)$,
\be
\label{eq4.19}
\nu_t = \exp (it\partial) \, .
\ee
Thus $\tau$ fails to be a trace because $\tau_0$ does. However we can 
compute the Hochschild coboundary $b \, \tau_0$, $b \, \tau_0 (a_0 , a_1) = 
\tau_0 (a_0 \, a_1) - \tau_0 (a_1 \, a_0)$. It only depends upon the symbols 
$\s (a_j) \in C^{\ify} (S^1)$ and is given by,
\be
\label{eq4.20}
b \, \tau_0 (a_0 , a_1) = \frac{-1}{2\pi i} \int \s (a_0) \, {\rm d} \s 
(a_1) \, .
\ee
One has 
$$
b \, \tau (a_0 , a_1) = \frac{1}{2\pi} \int b \, \tau_0 (a_0 (t) , a_1 (t)) 
\, {\rm d} t
$$
where $a(t) = \nu(t)(a)$, and for homogeneous elements of $\Ac$, $b \, \tau 
(a_0 , a_1) = 0$ unless the total degree is $(0,0)$. For such elements we 
thus get,
$$
b \, \tau (a_0 , a_1) = \frac{1}{2\pi} \int b \, \tau_0 (a_0 (t) , a_1 (t)) 
\, {\rm d} t = b \, \tau_0 (a_0 , a_1) = -\frac{1}{2\pi i} \int \s (a_0) \, 
{\rm d} \s (a_1)
$$
so that,
\be
\label{eq4.21}
b \, \tau (a_0 , a_1) = \frac{-1}{2 \pi i} \int \s (a_0) \, {\rm d} \s 
(a_1) \qquad \fl a_j \in \Ac \, .
\ee
Thus, even though $\tau$ is not a trace we do control by how much it fails 
to be a trace and this allows us to define a {\it cycle} in the sense of 
\cite{Coober2} using both first and second derivatives to define the 
differential,
\be
\label{eq4.22}
\Ac \build\longra_{}^{d} \Om^1 \, .
\ee
More precisely let us define the $\Ac$-bimodule $\Om^1$ with underlying 
linear space the direct sum, $\Om^1 = \Ac \op \Om^{(2)} (S^1)$ where 
$\Om^{(2)} (S^1)$ is the space of differential forms $f(\t) \, {\rm d} 
\t^2$ of weight 2 on $S^1$. The bimodule structure is defined by,
\begin{eqnarray}
\label{eq4.23}
a (\xi , f) = (a\xi , \s (a) f)\\ (\xi,f) b = (\xi b , -i \, \s (\xi) \, \s 
(b)' + f \s (b))
\nonumber 
\end{eqnarray}
for $a,b \in \Ac$, $\xi \in \Ac$ and $f \in \Om^{(2)} (S^1)$.

\smallskip

\noindent The differential $d$ of (\ref{eq4.22}) is then given by,
\be
\label{eq4.24}
da = \partial a + \frac{1}{2} \, \s (a)'' \, {\rm d} \t^2
\ee
as in a Taylor expansion.

\smallskip

\noindent The functional $\int$ is defined by,

\be
\label{eq4.25}
\int (\xi , f) = \tau (\xi) + \frac{1}{2\pi i} \int f \, {\rm d} \t \, .
\ee
We then have,

\medskip

\noindent {\bf Proposition 5.} {\it The triple $(\Om , d , \int)$ is a {\rm 
cycle}, i.e. $\Om = \Ac \op \Om^1$ equipped with $d$ is a graded 
differential algebra $($with $\Om^0 = \Ac)$ and the functional $\int$ is a 
closed graded trace on $\Om$.} 
\medskip

\noindent {\bf Proof.} One checks directly that $\Om^1$ is an 
$\Ac$-bimodule so that $\Om = \Ac \op \Om^1$ is a graded algebra. The 
equality $\s (\partial a) = i \, \s (a)'$ together with (\ref{eq4.23}) show 
that $d(ab) = (da) \, b + a \, db \quad \fl a,b \in \Ac$. It is clear also 
that $\int {\rm d} a = 0 \quad \fl a \in \Ac$. It remains to show that 
$\int$ is a (graded) trace, i.e. that $\int a \om = \int \om a \quad \fl 
\om \in \Om^1$, $a \in \Ac$.

\smallskip

\noindent With $\om = (a_1 , f)$ one has 
$$
a\om - \om a = (a a_1 - a_1 a , \s (a) f + i \s (a_1) \, \s (a)' - f \s 
(a)) = (a a_1 - a_1 a , i \s (a_1) \, \s (a)') \, .
$$
Thus (\ref{eq4.21}) shows that $\tau (a a_1 - a_1 a) + \frac{1}{2\pi} \int 
i\s (a_1) \, \s (a)' \, {\rm d} \t = 0$. $\Box$

\bigskip

\noindent We let $\chi$ be the cyclic 1-cocycle which is the character of 
the above cycle, explicitly,
\be
\label{eq4.26}
\chi (a_0 , a_1) = \int a_0 \, {\rm d} a_1 \qquad \fl a_0 , a_1 \in \Ac 
\, .
\ee
As above in proposition 3, we let $\vp$ be the cochain,
\begin{eqnarray}
\label{prop3}
&&\vp_0 (a) =  {\rm Trace}(\,a\, \vert D \vert^{-s})_{s=0},\nonumber \\
&& \vp_2(a_0,a_1,a_2)=\,\frac{1}{24}\int\!\!\!\!\!\!- \ a_0  \,\d(a_1) \,\d^2( a_2) \,  \, \vert 
D \vert^{-3} \,.\nonumber
\end{eqnarray}
\bigskip

\noindent {\bf Theorem 5.} {\it The local index formula of the spectral 
triple $(\Ac , \Hc , D)$ is given by the cyclic cocycle $\chi$ up to the coboundary of the cochain ($\vp_0$,$\vp_2$}).

\bigskip
\noindent In other words the cocycle $\psi_1$ of proposition 3
is equal to $\chi$. This follows from
(\ref{eq4.15}), (\ref{eq4.16}).

\noindent We leave it as an exercice for the reader
to compute the (non-zero) pairing between the above  cyclic cocycle $\chi$ 
and the K-theory class of the basic unitary,
\be
\label{eq3.22}
U = \left[ \matrix{\a &-q \b^* \cr \b &\a^* \cr} \right]
\ee
Applying (\ref{eta3})
 we
obtain the following corollary,
\bigskip

\noindent {\bf Corollary 1.} {\it The character ${\rm Trace}(a_0 [F, a_1])$ of the spectral 
triple $(\Ac , \Hc , D)$ is given by the cyclic cocycle $\chi$ up to the coboundary of the cochain $\psi_0$ given by}  $\psi_0 (a) =  2\,\,{\rm Trace}(\,a\,P \, \vert D \vert^{-s})_{s=0}$.
\bigskip

\noindent The cochain $\psi_0$ is only non-zero on elements which are functions
of $\b^* \b$ as one sees for homogeneity reasons using the bigrading and the natural
basis of $\Ac$ given by the $\a^k (\b^*)^n \b^m$. It is thus 
entirely determined by the values $\psi_0((\b^* \b)^n)$.
It is an interesting
problem to compute these functions of $q$. In order to state the result
we recall that the Dedekind eta-function for the modulus $q^2$
is given by,
\be
\label{eqded}
\eta(q^2)=\, q^{\frac{1}{12}} \prod_1^{\ify} (1-q^{2n}).
\ee
We let $G$ be its logarithmic derivative $q^2 \partial_{q^2} \eta(q^2)$, (up to sign and after substraction 
of the constant term),
\be
\label{eqgg}
G(q^2)=\,  \sum_1^{\ify}n\, q^{2n} (1-q^{2n})^{-1}.
\ee
\bigskip

\noindent {\bf Theorem 6.} {\it The functions $\frac{1}{2}\psi_0((\b^* \b)^r)$
of the variable $q$ are of the form $q^{-2r}(q^2\,R_r(q^2)-G(q^2))$ where 
$R_r(q^2)$ are rational fractions of $q^2$ with poles  only at roots of unity.
}

\bigskip
\noindent {\bf Proof.} The first step is to prove that the diagonal
terms $d(n,i,j)$ of the matrix $(\b^* \b)^r$ fulfill the equality,
\be
\label{diag}
d(n,n,j)=\,  \prod_{l=0}^{r-1}\, \frac{q^{2n+2j}- q^{4n+2+2l}}{1-q^{4n+4+2l}}.
\ee
This follows by writing $(\b^* \b)^r=\b^{*r} \b^r=(\b_+^*+\b_-^*)^r
(\b_++\b_-)^r$ and observing that, since $i=n$, the only term which contributes is 
$(\b_+^*)^r
(\b_+)^r$.

\noindent We then change variables as above replacing $n$ by $n/2$
and $j$ by $-n/2 +y$, which gives for the value of $\frac{1}{2}\psi_0((\b^* \b)^r)$
the value at $s=0$ of the sum,
\be
\label{reg}
Z(s)=\sum_{n=0}^{\ify}\,n^{-s} \sum_{y=0}^{n} \,  \prod_{l=0}^{r-1}\, \frac{q^{2y}- q^{2n+2+2l}}{1-q^{2n+4+2l}}.
\ee
(with the usual convention for $n=0$).

\noindent To understand the appearance of $G(q^{2})$ and the 
corresponding coefficient, let us take the constant term
in 
\be
\label{pol}
P(q^{2y})=\prod_{l=0}^{r-1}\, \frac{q^{2y}- q^{2n+2+2l}}{1-q^{2n+4+2l}}
\ee
 viewed as a polynomial in $q^{2y}$, which gives, with $x=q^{2n}$,
\be
\label{fra1}
R(x)=\prod_{l=0}^{r-1}\, \frac{(- q^{2+2l}x)}{1-q^{4+2l}x}
\ee
The fraction $R(1/z)$ has simple distinct poles at $z= q^{4+2l}$
and vanishes at $\ify$, thus we can express it in the 
form,
\be
\label{fra2}
R(1/z)=\sum_{l=0}^{r-1}\, \frac{\lb_l}{z-q^{4+2l}}
\ee
Each term  contributes to (\ref{reg})
by the value at $s= 0$ of,
\be
\label{reg1}
 \lb_l \sum_{n=0}^{\ify}\,n^{-s} (n+1) \,  \frac{q^{2n}}{1-q^{2n+4+2l}}.
\ee
The value at $s=0$ makes sense as a convergent series,
and the coefficient of $G(q^{2})$ is obtained by setting 
$n'=n+l+2$ which gives $\lb_l \, q^{-4-2l}\, G(q^{2})$.
Thus the overall coefficient for $G(q^{2})$ is, using (\ref{fra1}), (\ref{fra2}) and the behaviour for $z=0$,
\be
\label{reg1}
\sum_{l=0}^{r-1}\, \frac{\lb_l}{q^{4+2l}} =-q^{-2r}.  
\ee

\noindent One has $n+1=n'-l-1$ and the above terms also
generate a non-zero multiple of the function,
\be
\label{reg3}
 G_0(q^{2}) = \,  \frac{q^{2n}}{1-q^{2n}}.
\ee
The coefficient is given by,
\be
\label{reg4}
c_0= -\sum_{l=0}^{r-1}\, \lb_l  \,(l+1)   \, q^{-2l-4} .
\ee
We need to show that the other terms coming from the non-constant terms in $P(q^{2y})$ exactly cancell the above
multiple of $G_0(q^{2})$, modulo rational functions of $q^{2}$.

\noindent Using the $q^2$-binomial coefficients $(\begin{array}{c} r  \\  k  \\ \end{array})_{q^2}
$
, one obtains, with $P$ as in (\ref{pol}), that,
\be
\label{pol2}
P(z)= N(z)\,\prod_{l=0}^{r-1}\, (1-q^{2n+4+2l})^{-1}
\ee
where,
\be
\label{pol2}
N(z)= \,\sum_{k=0}^{r}\, (-1)^k q^{k(k+1)}q^{2k n}\, (\begin{array}{c} r  \\  k  \\ \end{array})_{q^2}
\, z^{r-k}
\ee
The constant term (in $z^0$) has already been taken care of, and for the others the effect of the summation $\sum_{y=0}^{n}$
in (\ref{reg}) is to replace $z^{r-k}$ in the above sum by,
\be
\label{pol3}
\sum_{y=0}^{n}\, q^{2y(r-k)}= \, \frac{1-q^{2(r-k)(n+1)}}{1-q^{2(r-k)}}
\ee
Thus the contribution of the other terms is governed by the rational fraction of $x=q^{2n}$
\be
\label{pol2}
Q(x)= \,(\sum_{k=0}^{r-1}\, (-1)^k q^{k(k+1)}\, (\begin{array}{c} r  \\  k  \\ \end{array})_{q^2}
 \frac{x^{k }-x^{r}q^{2(r-k)}}{1-q^{2(r-k)}}
)\,\prod_{l=0}^{r-1}\, (1-x \,\,q^{4+2l})^{-1}
\ee
The degree of the numerator is the same as the degree of the
denominator, all poles are simple, and we can thus expand
$Q(x)$ as,
\be
\label{fra4}
Q(x)=\, \mu + \sum_{l=0}^{r-1}\, \frac{\mu_l}{1-x \, q^{4+2l}}
\ee
The same reasoning as above shows that modulo rational functions
of $q^2$, each term contributes to (\ref{reg}) by a multiple 
of $G_0(q^{2})$, while the overall coefficient is the sum of the $\mu_l$,
\be
\label{coeff}
c_1=Q(0)-\, Q(\ify)
\ee
Using (\ref{pol2}) one obtains,
\be
\label{pol4}
c_1= \,(1-q^{2r})^{-1}+(-1)^r q^{-r(r+1)}\sum_{k=0}^{r-1}\, (-1)^k \,(\begin{array}{c} r  \\  k  \\ \end{array})_{q^2}
 \, \frac{q^{k(k-1)}}{1-q^{2(r-k)}}
\ee
To compute the $\lb_l$ one takes the residues of (\ref{fra2})
which gives the formula,
\be
\label{pol5}
\lb_j= -\,(-1)^j \, q^{(4+j^2+j(3-2r)-3r+r^2)} \,
(\begin{array}{c} r-1\\  j  \\ \end{array})_{q^2}
\rho(r-1)^{-1}
\ee
where,
\be
\label{fac}
\rho(r-1)= \prod_{a=1}^{r-1}\,(q^{2a} -1)
\ee
This gives the following formula for the coefficient $c_0$,
\be
\label{c0}
c_0=\rho(r-1)^{-1} \sum_{j=0}^{r-1}\,(-1)^j \,(j+1)\,\,
q^{(j^2+j-2rj-3r+r^2)}\, (\begin{array}{c} r-1\\  j  \\ \end{array})_{q^2}
\ee
The fundamental cancellation now is the
identity 
\be
\label{cancell}
c_0 +c_1 =0
\ee
which
is proved by differentiation of the 
$q^2$-binomial formula.

\noindent
The above discussion provides an explicit formula for the rational fractions $R_r$ which allows to check that their
only poles are roots of unity.
 $\Box$

\noindent The simple expression $q^{-2r}G(q^2)$ blows up exponentially for $r \ra \ify$ and if it were alone it would be impossible to extend the cochain $\psi_0$ from the 
purely algebraic to the smooth framework. However, 
$$
\psi_0((\b^* \b)^r) \ra 1+ 2q^2/(q^2 -1)
\,\,{\rm when }\,\, r \ra \ify.
$$
Thus it is only by the virtue of the rational approximations $q^2\,R_r(q^2)$
of $G(q^2)$ that the tempered behaviour of $\psi_0((\b^* \b)^r)$ is insured.

\noindent The list of the first $R_r(q)$ is as follows,
\begin{eqnarray}
\label{ratfracs}
&&R_1[q]=\frac{3}{2(1-q)},\qquad
R_2[q]=\frac{2+5q-3q^2}{2(-1+q)^2(1+q)}, \\
&&R_3[q]=\frac{2+8q+13q^2+11 q^3 -q^4 -3q^5}{2(-1+q^2)^2(1+q+q^2)},\nonumber \\
&&R_4[q]=\frac{2+10q+24q^2+43 q^3 +50q^4 +46q^5+24q^6+4q^7-4q^8-3q^9}{2(1+q^2)(-1-q+q^3+q^4)^2}.\nonumber
\end{eqnarray}

\noindent Finally, note that the appearance of the function 
$G(q^2)$ in $\psi_0$ is not an artefact which could be eliminated by a better choice of 
 cochain with the same coboundary. Indeed since 
$$
\a \a^* -\a^* \a = \, (1-q^2) \, \b^* \b
$$
the coboundary $b\psi_0(\a ,\a^*)$ still involves $G(q^2)$.

\bigskip

\section{Quantum groups and invariant cyclic cohomo\-logy}

The main virtue of the above spectral triple for ${\rm SU}_q (2)$ is its 
invariance under left translations (cf. \cite{CP}). More precisely the following equalities 
define an action of the envelopping algebra $\Uc = U_q ({\rm SL} (2))$ on 
$\Hc$ which commutes with $D$ and implements the translations on $C^{\ify} 
({\rm SU}_q (2))$,
\be
\label{eq5.1}
k \, e_{ij}^{(n)} = q^j \, e_{ij}^{(n)}
\ee
\be
\label{eq5.2}
e \, e_{ij}^{(n)} = q^{-n+\frac{1}{2}} (1-q^{2(n+j+1)})^{1/2} 
(1-q^{2(n-j)})^{1/2} (1-q^2)^{-1} \, e_{ij+1}^{(n)}
\ee
while $f=e^*$.

\smallskip

\noindent With these notations one has,
\be
\label{eq5.3}
ke = q\,e\,k \, , \ kf = q^{-1} fk \, , \ [e,f] = \frac{k^2 - k^{-2}}{q-q^{-1}} 
\, .
\ee
The vector $\Om= \,e_{(0,0)}^{(0)}$ is preserved by the action and one has a 
natural densely defined action of $\Uc$ on $\Ac = C^{\ify} ({\rm SU}_q (2))$ such that,
\be
\label{eq5.4}
h(x) \Om = h(x\Om) \qquad \fl x \in \Ac \, , \ h \in \Uc \, .
\ee
The coproduct is given by,
\be
\label{eq5.5}
\D k = k \ot k \, , \ \D e = k^{-1} \ot e + e \ot k \, , \ \D f = k^{-1} \ot 
f + f \ot k
\ee
and the action of $\Uc$ on $\Ac$ fulfills,
\be
\label{eq5.6}
h(xy) = \sum h_{(1)} (x) \, h_{(2)} (y) \qquad \fl x,y \in \Ac \, ,
\ee
$\fl h \in \Uc$ with $\D h = \sum h_{(1)} \ot h_{(2)}$.

\smallskip

\noindent On the generators $\a , \a^* , \b , \b^*$ of $\Ac$ one has,
\be
\label{eq5.7}
k(\a) =  q^{-1/2} \a \, , \ k(\b) = q^{-1/2} \b \, , \ e(\a) = q\b^* \, ,  
\ee
$$
e(\b) = -\a \, , \ e(\a^*) = 0 \, , \ e(\b^*) = 0 \, .
$$
This representation of $\Uc$ in $\Hc$ generates the regular representation 
of the compact quantum group ${\rm SU}_q (2)$ and we let $M = \Uc''$ be the 
von~Neumann algebra it generates in $\Hc$. It is a product $M = 
\build\prod_{n \in \frac{1}{2} \Nb}^{} M_{2n+1} (\Cb)$ of matrix algebras, 
where $M_{2n+1} (\Cb)$ acts with multiplicity $2n+1$ in the space $\Hc_{(n)} 
= \{\hbox{span of} \ e_{i,j}^{(n)}\}$.

\smallskip

\noindent The elements of $\Uc$ are unbounded operators affiliated to $M$ 
and at the qualitative level we shall leave the freedom to choose a weakly 
dense subalgebra $\Cc$ of $M$. Since all the constructions performed so far 
in this paper were canonically dependent on the spectral triple, the ${\rm 
SU}_q (2)$ equivariance of $(\Ac , \Hc , D)$ should entail a corresponding 
{\it invariance} of all the objects we delt with. We shall concentrate on 
the cyclic cohomology aspect and show that indeed there is a fairly natural 
and simple notion of {\it invariance} fulfilled by all cochains involved in 
the above computation.

\smallskip

\noindent The main point is that we can enlarge the algebra $\Ac$ to the 
algebra $\Dc = \Ac \semi 
\Cc$ generated by $\Ac$ and $\Cc$, extend the cochains on $\Dc$ by similar formulas and use the commutation,
\be
\label{eq5.8}
[D,c] = 0 \qquad \fl c \in \Cc
\ee
to conclude that the extended cochains fulfill the following key property,

\medskip

\noindent {\bf Definition 2.} {\it Let $\Dc$ be a unital algebra, $\Cc \sbs 
\Dc$ a (unital) subalgebra and $\vp \in C^n (\Dc)$ an} n-{\it cochain. We shall 
say that $\vp$ is {\rm $\Cc$-constant} iff both $\vp (a^0 , \ldots , a^n)$ 
and $(b\vp) (a^0 , \ldots , a^{n+1})$ vanish if one of the $a^j, \, \,j \geq 
1$ is in $\Cc$.}

\medskip

\noindent When $\Cc = \Cb$ this is a normalization condition.

\smallskip

\noindent When $\vp$ is $\Cc$-constant then $B_0 \vp (a^0 , \ldots , 
a^{n-1}) = \vp (1 , a^0 , \ldots , a^{n-1})$ so that $B \vp = A B_0 \vp$ is 
also $\Cc$-constant. It follows that $\Cc$-constant cochains form a 
subcomplex of the $(b,B)$ bicomplex of $\Dc$ and we can develop cyclic 
cohomology in that context, parallel to (\cite{Coober2},\cite{Coober3}). We shall denote by 
$HC_{\Cc}^* (\Dc)$ the corresponding theory. 

\noindent In the above context we take for $\Dc$ the algebra $\Ac \semi \Cc$
and use the lighter notation $HC_{\Cc}^* (\Ac)$ for the corresponding theory. 
 
\noindent Let us now give examples of specific cochains on $\Ac = C^{\ify} ({\rm SU}_q 
(2))$ which extend to $\Cc$-constant cochains on $\Ac \semi \Cc = \Dc$. We 
start with the non local form of the Chern character of the spectral triple,
\be
\label{eq5.9}
\psi_1 (a^0 , a^1) = {\rm Trace} \, (a^0 \, [F,a^1]) \qquad \fl a^0 , a^1 
\in \Ac \, .
\ee
Let us show how to extend $\psi_1$ to an $M$-constant cochain on $\Ac \semi M 
= \Dc$. An element of $\Dc$ is a finite linear combination of monomials $a_1 
\, m_1 \, a_2 \, m_2 \ldots a_{\ell} \, m_{\ell}$, where $a_j \in \Ac$, 
$m_{\ell} \in M$. But $[F,a]$ is a trace class operator for any $a \in \Ac$, 
while $[F,m] = 0 \quad \fl m \in M$, thus we get,
\be
\label{eq5.10}
[F,x] \in \Lc^1 \qquad \fl x \in \Dc = \Ac \semi M \, .
\ee
We can thus define $\wt{\psi}_1$ as the character of the module $(\Hc , F)$ 
on $\Dc$, namely,
\be
\label{eq5.11}
\wt{\psi}_1 (x_0 , x_1) = {\rm Trace} \, (x_0 \, [F , x_1]) \, .
\ee
It is clear that $\wt{\psi}_1$ is $M$-constant and that $b \wt{\psi}_1 = 0$ 
so that $b \wt{\psi}_1$ is also $M$-constant, $\wt{\psi}_1 \in HC_M^1 
(\Ac)$. This example is quite striking in that we could extend $\psi_1$ to a 
very large algebra. Indeed if we stick to bounded operators $M$ is the 
largest possible choice for $\Cc$. A similar surprising extension of a
cyclic $1$-cocycle in a von-Neumann algebra context already
occured in the anabelian $1$-traces of (\cite{trans}).
\noindent As a next example let us take the functional on $\Ac$ which is the natural 
trace,
\be
\label{eq5.12}
\psi_0 (x) = \frac{1}{2\pi} \int \s (x) \, {\rm d} \t \qquad \fl x \in 
C^{\ify} ({\rm SU}_q (2)) \, .
\ee
When written like this, its ${\rm SU}_q (2)$-invariance is not clear and in 
fact cannot hold in the simplest sense since this would contradict the 
uniqueness of the Haar state on $C^{\ify} ({\rm SU}_q (2))$. Let us however 
show that $\psi_0$ extends to an $M$-constant cochain (in fact an 
$M$-constant trace) on $\Dc = \Ac \semi M$ as above. To do this we rewrite 
(\ref{eq5.12}) as,
\be
\label{eq5.13}
\psi_0 (x) = {\rm Tr}_{\om} (x \, \vert D \vert^{-3}) \qquad \fl x \in 
C^{\ify} ({\rm SU}_q (2))
\ee
where ${\rm Tr}_{\om}$ is the Dixmier trace (\cite{dix})(\cite{book}) and simply write the 
extension as,
\be
\label{eq5.14}
\wt{\psi}_0 (x) = {\rm Tr}_{\om} (x \, \vert D \vert^{-3}) \, .
\ee
For any monomial $\mu = a_1 \, m_1 \ldots a_m \, m_m$ as above one has 
$[D,\mu]$ bounded and $[\vert D \vert , \mu]$ bounded. Thus it follows from 
the general properties of ${\rm Tr}_{\om}$ that,
\be
\label{eq5.15}
\wt{\psi}_0 (xy) = \wt{\psi}_0 (yx) \qquad \fl x,y \in \Dc = \Ac \semi M \, .
\ee
This shows of course that $\wt{\psi}_0$ is a 0-cycle in the invariant cyclic 
cohomology $HC_M^0 (\Ac)$. After giving these simple examples it is natural 
to wonder wether the above notion of $\Cc$-constant cochain is restrictive 
enough. Here is a simple consequence of this hypothesis:

\medskip

\noindent {\bf Proposition 6.} {\it Let $\Cc \sbs \Dc$ be unital algebra and 
$\vp \in C_{\Cc}^n$ be a $\Cc$-constant cochain on $\Dc$. Then for any 
invertible element $u \in \Cc$ one has,
$$
\vp (u \, a^0 u^{-1} , u \, a^1 u^{-1} , \ldots , u \,  a^n u^{-1}) = \vp (a^0 , 
\ldots , a^n) \, .
$$
}

\noindent {\bf Proof.} One has $b \vp (a^0 , u , a^1 , \ldots , a^n) = 0$ so 
that $\vp (a^0 u , a^1 , \ldots , a^n) - \vp (a^0 , u a^1 ,$ $\ldots , a^n) 
= 0$ since all other terms have $u$ as an argument and hence vanish. 
Similarly $\vp (a^0 , \ldots , a^{j-1} u , a^j , \ldots ,a^n) = \vp (a^0 , 
\ldots , a^{j-1}, ua^j , \ldots , a^n)$ for all $j \in \{ 1,\ldots , n\}$ 
and $\vp (u a^0 , \ldots , a^n) = \vp (a^0 , \ldots , a^n u)$. Applying 
these equalities yields the statement. $\Box$

\bigskip

\noindent Let us now consider the more sophisticated cochains which appeared 
throughout and show how to extend them to $\Cc$-constant cochains on $\Dc = 
\Ac \semi \, \Cc$ for suitable algebra $\Cc$ describing the quantum group 
${\rm SU}_q (2)$.

\smallskip

\noindent We first note that the action of the envelopping algebra $\Uc = 
U_q ({\rm SL} (2))$ on $\Ac$ extends to an action on the algebra of 
pseudo-differential operators. First it extends to $\Bc$ with the following 
action on the generators $\a_{\pm}$, $\a_{\pm}^*$, $\b_{\pm}$, $\b_{\pm}^*$,
\be
\label{eq5.16}
k (\a_{\pm}) = q^{-1/2} \a_{\pm} \, , \ k(\b_{\pm}) = q^{-1/2} \b_{\pm} \, ,
\ee
and
\be
\label{eq5.17}
e(\a_{\pm}) = q \b_{\mp}^* \, , \ e(\b_{\pm}) = - \a_{\mp} \, , \ e 
(\a_{\pm}^*) = 0 \, , \ e (\b_{\pm}^* ) = 0 \, .
\ee
Moreover $\Uc$ acts through the trivial representation on $D$, $\vert D 
\vert$ and $F$.

\smallskip

\noindent In fact it is important to describe the action of $\Uc$ on 
arbitrary pseudo-differential operators by a closed formula and this is 
achieved by,

\medskip

\noindent {\bf Proposition 7.} {\it The action of the generators $k,e,f$ of 
$\Uc$ on pseudo-differential operators $P$ is given by,} a) $k(P) = 
kPk^{-1}$, b) $e(P) = e P k^{-1} - q k^{-1} P e$, \break c) $f(P) = f P k^{-1} 
- 
q^{-1} k^{-1} P f$.

\medskip

\noindent {\bf Proof.} These formulas just describe the tensor product of 
the action of $\Uc$ in $\Hc$ by the contragredient representation, since the 
antipode $S$ in $\Uc$ fulfills
\be
\label{eq5.18}
S(k) = k^{-1} \, , \ S(e) = -qe \, , \ S(f) = -q^{-1} f \, .
\ee
One checks directly that they agree with (\ref{eq5.16}) and (\ref{eq5.17}) 
on the generators $\a_{\pm} , \ldots , \b_{\pm}^*$ as well as on $D$, $\vert 
D \vert$ and $F$. Thus we are just using the natural implementation of the 
action of $\Uc$ which extends this action to operators. $\Box$

\bigskip

\noindent The only technical difficulty is that the generators of $\Uc$ are 
unbounded operators in $\Hc$ so that to extend cochains to $\Ac \semi \, 
\Uc$ requires a little more work. In fact the only needed extension is for 
the residue,
\be
\label{eq5.19}
\int\!\!\!\!\!\!- \, P = {\rm Res}_{s=0} {\rm Trace} \, (P \, \vert D 
\vert^{-s}) \, .
\ee
Using formula $(\t)$ we can reexpress (\ref{eq5.19}) as follows,
\be
\label{eq5.20}
\int\!\!\!\!\!\!- \, P = \frac{1}{2} \ \hbox{coefficient of} \ \log t^{-1} \ 
\hbox{in} \ {\rm Trace} \, (P e^{-t D^2}) \, .
\ee
Thus more precisely we let $\t_P (t) = {\rm Trace} \, (P e^{-t D^2})$ and 
assume that it has an asymptotic expansion for $t \ra 0$ of the form 
\be
\label{eq5.21}
\t_P (t) \sim \sum a_{\a} \, t^{-\a} + \lb \log t^{-1} + a_0 + \cdots \, ,
\ee
then the equality between (\ref{eq5.19}) and (\ref{eq5.20}) holds, both 
formulas giving $\lb / 2$. In our context we could use (\ref{eq5.20}) above 
instead of (\ref{eq5.19}) since we always controlled the size of $\Zc_b (s)$ 
on vertical strips to perform the inverse Mellin transform.

\smallskip

\noindent Let now $L$ be an arbitrary extension of the linear form on 
function $f \in C^{\ify}$ $(] 0,\ify[)$ which satisfies,
\be
\label{eq5.22}
L(f) = \frac{1}{2} \ \hbox{coefficient of} \ \log t^{-1} \ \hbox{if} \ f \ 
\hbox{admits}
\ee
$$
\hbox{an asymptotic expansion (\ref{eq5.21})}.
$$
We then extend the definition (\ref{eq5.19}) by,
\be
\label{eq5.23}
{\int\!\!\!\!\!\!-}_L \, P = L (\t_P (t)) \, .
\ee
With these notations we then have,

\medskip

\noindent {\bf Proposition 8.} {\it Let $(k_1 , \ldots , k_n)$ be a 
multi-index, then the formula 
$$
\wt\psi (a^0 , \ldots , a^n) = {\int\!\!\!\!\!\!-}_L \, a^0 [D,a^1]^{(k_1)} 
\ldots [D,a^n]^{(k_n)} \vert D \vert ^{-n-\vert k \vert}
$$
where $T^{(k)} = \d^k (T)$, defines a $\Uc$-constant extension of the 
restriction $\psi$ of $\wt\psi$ to $\Ac$ to the algebra $\Dc = \Ac \semi \, 
\Uc$.}

\medskip

\noindent {\bf Proof.} In computing $b \wt\psi$ one uses the equality
\begin{eqnarray}
\label{eq5.24}
\d^k ([D,a \,b]) &= &\d^k ([D,a]) \, b + a \, \d^k ([D,b]) + \sum_{j=0}^{k-1} C_k^j \, 
\d^j ([D,a]) \, \d^{k-j} (b) \nonumber \\
&+ &\sum_{\ell = 1}^k \d^{\ell} (a) \, \d^{k-\ell} ([D,b]) \, . 
\end{eqnarray}
Thus in $b \, \wt\psi (a_0 , \ldots , a_{n+1})$ the only term which does not 
involve a derivative of $a$ is of the form,
\be
\label{eq5.25}
\int\!\!\!\!\!\!- \, a_{n+1} \, T \vert D \vert^{-n-\vert k \vert} - 
\int\!\!\!\!\!\!- \, T \vert D \vert^{-n-\vert k \vert} \, a_{n+1} \, .
\ee
This shows that $b \, \wt\psi$ vanishes if any of the $a_j \in \Uc$ for $j = 1 
, \ldots , n$. For $j = n+1$, i.e. for $a_{n+1} = v \in \Uc$ one has the 
term (\ref{eq5.25}) but since $v$ commutes with $D$ one has,
\be
\label{eq5.26}
\t_{vT} = \t_{Tv} \, ,
\ee
and one gets the desired result. $\Box$

\bigskip

\noindent This proposition shows the richness of the space of $\Uc$-constant 
cochains, but it does not address the more delicate issue 
of computing the cyclic cohomology $HC_{\Uc}^*$ $(\Ac)$. A much more careful choice of $L$ would be 
necessary if one wanted to lift cocycles to cocycles. 

\smallskip

\noindent We shall now show that $HC_{\Uc}^* (\Ac)$ which obviously maps to 
 the ordinary cyclic theory $HC^* (\Ac)$,
\be
\label{eq5.27}
HC_{\Uc}^* (\Ac) \build\longra_{}^{\rho} HC^* (\Ac) \, ,
\ee
 also maps in fact to the ``twisted'' cyclic cohomology 
$HC_{\rm inv}^* (\Ac, \, \theta)$ proposed in 
\cite{Tu}, where $\theta$ is the inner automorphism
implemented by $k^2$. This will allow to put the latter proposal in the correct 
perspective. Indeed the drawback of this simple variation on (\cite{Coober3})
is that it lacks the relation to $K$-theory which is the back-bone
of cyclic cohomology. This was a good reason to refrain from developping
such a "twisted" form of the general theory
in spite of its previous appearance in (\cite{CM1} cf. equation 2.28 p.14)
and of its merit which is to connect with the various "differential calculi" on quantum groups (\cite{woro},\cite{woro2}). However the next proposition shows
that 
it would be very interesting to use it as a "detector" of classes in $HC_{\Uc}^* (\Ac)$.

\smallskip

\noindent To see what happens, let us start with a $\Uc$-constant 
$0$-dimensional cochain $\psi$ on $\Ac \semi \, \Uc$ and get an analogue of 
the group invariance provided by proposition 6. One has of course $\psi 
(kak^{-1}) = \psi (a)$ but this is not much. We would like a similar 
statement for the other generator $e$ of $U_q(SL(2)$.
 Now by proposition 7 one has $e(a) = eak^{-1} - qk^{-1} ae$ 
so that $e(a)k^2 = eak - qk^{-1} aek^2$. But $\Uc$ is in the centraliser of 
$\psi$ by proposition 6 and thus,
$$
\psi (eak) = \psi (kea) \, , \ \psi (k^{-1} a e k^2) = \psi (eka) \, ,
$$
hence $\psi (e(a) \, k^2) = 0$. One gets in general,
\be
\label{eq5.28}
\psi (h(a)k^2) = \ve (h) \, \psi (a \,k^2) 
\ee
which is the usual invariance of a linear form. More generally one has,

\medskip

\noindent {\bf Proposition 9.}  {\it The equality $\rho_{\theta} (\psi)(a_0 , \ldots , 
a_n) = \psi (a_0 , \ldots , a_n k^2)$ defines a morphism,}
$$
HC_{\Uc}^* (\Ac) \build\longra_{}^{\rho_{\theta}} HC_{\rm inv}^* (\Ac, \, \theta) \, .
$$
 {\it where $\theta$ is the inner automorphism 
implemented by $k^2$. }
\medskip

\noindent We have seen above (\ref{eq5.11}
),(\ref{eq5.14}) that the basic cohomology classes in
the ordinary cyclic theory $HC^* (\Ac)$ of $\Ac$ lift to 
actual cocycles in $HC_{M}^* (\Ac)$ where $M$ is the von-Neumann
algebra bicommutant of $\Uc$. It is however not clear that they
lift to $HC_{\Uc}^* (\Ac)$ since the generators of $\Uc$ are
unbounded operators. We can however insure that such liftings 
exist in the entire cyclic cohomology (\cite{theta})(\cite{JLO})
since the $\theta$-summability of the spectral triple 
continues to hold for the algebra $\Ac \semi \Uc$. This point is not
unrelated to
the attempt by Goswami in (\cite{G1}).

\noindent What
 we have shown here is that the local formulas work perfectly well
in the context of quantum groups,
and that the framework of NCG needs
no change whatsoever, at least as far as $SU_q(2)$
is concerned. The only notion that requires more work is that of invariance
in the q-group context. 

\noindent  Finally the above notion of invariant cyclic cohomology
is complementary to the theory developped in (\cite{CM2},\cite{CM1}).
In the latter the Hopf action is used to construct ordinary cyclic cocycles
from twisted-traces.
In the q-group situation, cocycles thus constructed from the right
translations should be left-invariant in the above sense.

\section{Appendix}

\noindent We have not defined carefully the smooth algebras
$C^{\ify}$ involved in section 6. A careful definition can 
however be deduced from their structure and the exact sequence
involving   $OP^{-\ify} $ and the symbol map provided by lemma 1.
What really matters is that the obtained algebras are stable
under
holomorphic functional calculus (h.f.c.)
and we shall now provide the technical lemma which allows
to check this point.

\noindent Let   $(B,\Hc,D)$ be a spectral triple.
As above we say that an operator $P$ 
in $\Hc$ is of order $\a$ when,
\be
\label{eq3.1}
\vert D \vert^{-\a} \, P \in \bigcap_{n=1}^{\ify} {\rm Dom} \, \d^n
\ee
where $\d$ is the unbounded derivation given by,
\be
\label{eq3.2}
\d (T) = \vert D \vert \, T - T \, \vert D \vert \, .
\ee
Thus $OP^0 = \build\bigcap_{n=1}^{\ify} {\rm Dom} \, \d^n$ is the algebra 
of operators of order $0$ and $OP^{-\ify} $ is a two sided ideal in $OP^0$.

\noindent  Let $\rho: B \ra C$
a morphism of $C^*$-algebras, $\Cc \sbs C$ be a subalgebra stable
under h. f. c. and $\lb: \Cc \ra \Lc(\Hc)$ be a linear map such that  $\lb(1)=1$ and,

\begin{eqnarray}
\label{hfc}
&&\lb (c) \in  OP^0, \,\,\forall c \in 
\Cc
 \\
&& \lb (a \,b)- \lb (a)  \lb (b) \in  OP^{-\ify},
\,\,\,\,\forall a, b \in 
\Cc
 \,.\nonumber
\end{eqnarray}
We then have the following,

\medskip

\noindent {\bf Lemma 2.} {\it Let $\Bc =\{x \in B\,; \,x \in OP^0, \,\rho(x)
\in \Cc,\, x- \lb( \rho(x)) \in  OP^{-\ify}\}$. Then $\Bc \sbs B$ is a 
subalgebra stable under holomorphic functional calculus.}
\medskip

\noindent {\bf Proof.} Let $x \in \Bc $ be invertible in $B$, let us show that
$x^{-1} \in \Bc $. Let $a=\rho(x)$, then since $\Cc$ is stable under h.f.c.
the inverse $b=\rho(x^{-1})$ of $a$ belongs to $\Cc$. Also since $x \in OP^0$
we have $x^{-1} \in OP^0$. Let us show that $x^{-1}- \lb(b) \in OP^{-\ify}$.
Since $a\,b=1$ one has by (\ref{hfc}), $\lb(a)\lb(b)-1 \in OP^{-\ify}$. But $x-
\lb(a) \in OP^{-\ify}$ and $OP^{-\ify}$ is a two-sided
ideal in $OP^0$, thus multiplying $x-
\lb(a)$ by $\lb(b)$ on the right,
we get $x \,\lb(b)-1 \in OP^{-\ify}$. Finally since
$x^{-1} \in OP^0$ we get, multiplying $x \,\lb(b)-1 $ on the left by $x^{-1}$
that $\lb(b)-x^{-1} \in OP^{-\ify}$.

\end{document}
\footnote{Coll\`ege de France, 3, rue Ulm, 75005 PARIS\\
et \\ I.H.E.S., 35, route de Chartres, 91440 BURES-sur-YVETTE }
\affiliation{\affnote{}Coll\`ege de France, 3, rue Ulm, 75005 PARIS\\
et \\ I.H.E.S., 35, route de Chartres, 91440 BURES-sur-YVETTE }